\documentclass[12pt]{amsart}
\usepackage{amsthm,amsmath,amssymb}
\usepackage[colorlinks=true,citecolor=blue,linkcolor=black,urlcolor=blue]{hyperref}

\usepackage{geometry}
\usepackage{graphicx,cite,enumitem}
\usepackage{color}

\numberwithin{equation}{section}
\numberwithin{figure}{section}

\theoremstyle{plain} 
\newtheorem{thm}{Theorem}

\theoremstyle{remark}

\theoremstyle{plain} 
\newtheorem{prob}{Problem}

\newcommand{\M}{\operatorname{T}}
\newcommand{\MM}{\operatorname{M}}
\newcommand{\Hf}{\operatorname{\textbf{H}}}
\newcommand{\wt}{\operatorname{wt}}

\newcommand{\tr}{\operatorname{tr}}

\begin{document}

\title{Problems in the Enumeration of Tilings}

\author{Tri Lai}
\address{Department of Mathematics, University of Nebraska -- Lincoln, NE 68588, U.S.A.}
\email{tlai3@unl.edu}
\thanks{This research was supported in part  by Simons Collaboration Grant (\# 585923).}



\date{\today}

\dedicatory{}

\begin{abstract}
Enumeration of tilings is the mathematical study concerning the total number of coverings of regions by similar pieces without gaps or overlaps.  Enumeration of tilings has become a vibrant subfield of combinatorics with connections and applications to diverse mathematical areas.  In 1999, James Propp published his well-known list of 32 open problems in the field. The list has got much attention from experts around the world. After two decades, most of the problems on the list have been solved and generalized. In this paper, we propose a set of new tiling problems. This survey paper contributes to the \emph{Open Problems in Algebraic Combinatorics 2022} conference (OPAC 2022) at the University of Minnesota.
\end{abstract}

\maketitle
\section{Introduction}\label{Sec:Intro}

Enumeration of tilings is a subfield of combinatorics studying the total number of coverings (called ``\emph{tilings}") of regions by similar pieces without gaps or overlaps. The first major result of the enumeration of tilings is usually credited to Percy Alexander MacMahon (1854--1929) with his beautiful formula for the number of lozenge tilings of a hexagon in the triangular lattice. However, MacMahon did not study tilings; instead, he worked on the enumeration of plane partitions. More than 100 years ago, he proved his celebrated theorem on the number of plane partitions fitting in a given box \cite{Mac}. Much later (in the 1980s), G. David and C. Tomei showed a simple bijection between lozenge tilings of a centrally symmetric hexagon of side-lengths $a,b,c,a,b,c$ (in cyclic order) and plane partitions fitting in an $a\times b\times c $-box \cite{DT}, as in Figure \ref{Macfig}. (Strictly speaking, David and Tomei only showed the bijection for the case $a=b=c$; however, their bijection could be easily extended for the general case.) This way, MacMahon's theorem implies a product formula for the tiling number of a hexagon. Since  then, MacMahon has been considered as one of the founding fathers of the field.

Actually, tiling-counting problems in the square lattice have already been investigated here and there for decades before the David--Tomei bijection, say under the form of small mathematical puzzles. For instance, no one really knows the author of the folklore puzzle on the number of ways to cover a rectangular stripe of width $2$  by domino pieces\footnote{The answer is the Fibonacci number.}.  Several results of the same flavor appeared in recreational and discrete mathematics, for instance, the work of David Klarner in the 1960s \cite{Klarner}.

A couple of significant results in the enumeration of tilings come from statistical physics. In the early 1960s, the physicists P.W. Kastelyn \cite{Kes} and H.N.V. Temperley and M.E. Fisher \cite{Fish} independently found an explicit formula for the number of dimer configurations of a rectangular grid graph. This result equivalently gives the enumeration of domino tilings of a rectangle. It would be a flaw here if we do not mention the well-known result of Fisher and Stephenson \cite{FS} concerning interactions of holes in dimer systems on the square lattice. Similarly, we cannot ignore the 1990 classical articles of J.H. Conway and J.C. Lagarias \cite{Conway} and W. Thurston \cite{Thurston} investigating the connection between tilings and group theory.

In 1999, James Propp published his well-known article, ``\emph{Enumeration of Matchings: Problems and Progress}" \cite{Propp}, tracking the progress on a list of 32 open problems in the field. He presented this list in a 1996 lecture as part of the special program on algebraic combinatorics organized at MSRI.  In a review on MathSciNet of the American Mathematical Society (AMS), Christian Krattenthaler (Professor at the University of Vienna) wrote about this list of problems: ``\emph{This list of problems was very influential; it called forth tremendous activity, resulting in the solution of several of these problems (but by no means all), in the development of interesting new techniques, and, very often, in results that move beyond the problems.}'' The enumeration of tilings has become a vibrant subfield of enumerative and algebraic combinatorics with connections and applications to diverse areas of mathematics, including representation theory, linear algebra, cluster algebra, group theory, mathematical physics, graph theory, probability, and discrete dynamical systems, just to name a few.
We also refer the reader to the excellent survey (also by Propp) \cite{Propp2} for many connections and applications of the enumeration of tilings.

As most of Propp's problems have been solved, we would like to propose a new set of  tiling problems. This batch of problems is independent of Propp's list; we do not include here the still-open problems from the 1999 list. The author does not attempt to collect all open problems in the field of enumeration of tilings. The choice of problems in this paper reflects the author's personal taste.

\section{Weighted Enumerations of Lozenge Tilings}
\begin{figure}\centering
\includegraphics[width=8cm]{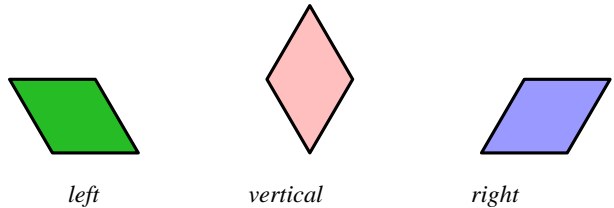}
\caption{Three orientations of the lozenges.}\label{Fig:orientation}
\end{figure}

Weighted enumeration is usually more challenging and often gives more insights than `plain' counting. This section is devoted to the weighted enumerations of lozenge tilings of regions in the (regular) triangular lattice. We orient the triangular lattice so that it accepts horizontal lattice lines. The \emph{lozenges} (the unions of two adjacent unit triangles) in the triangular lattice have three possible orientations: left, vertical, and right, as in Figure \ref{Fig:orientation}.

\begin{figure}\centering
\includegraphics[width=14cm]{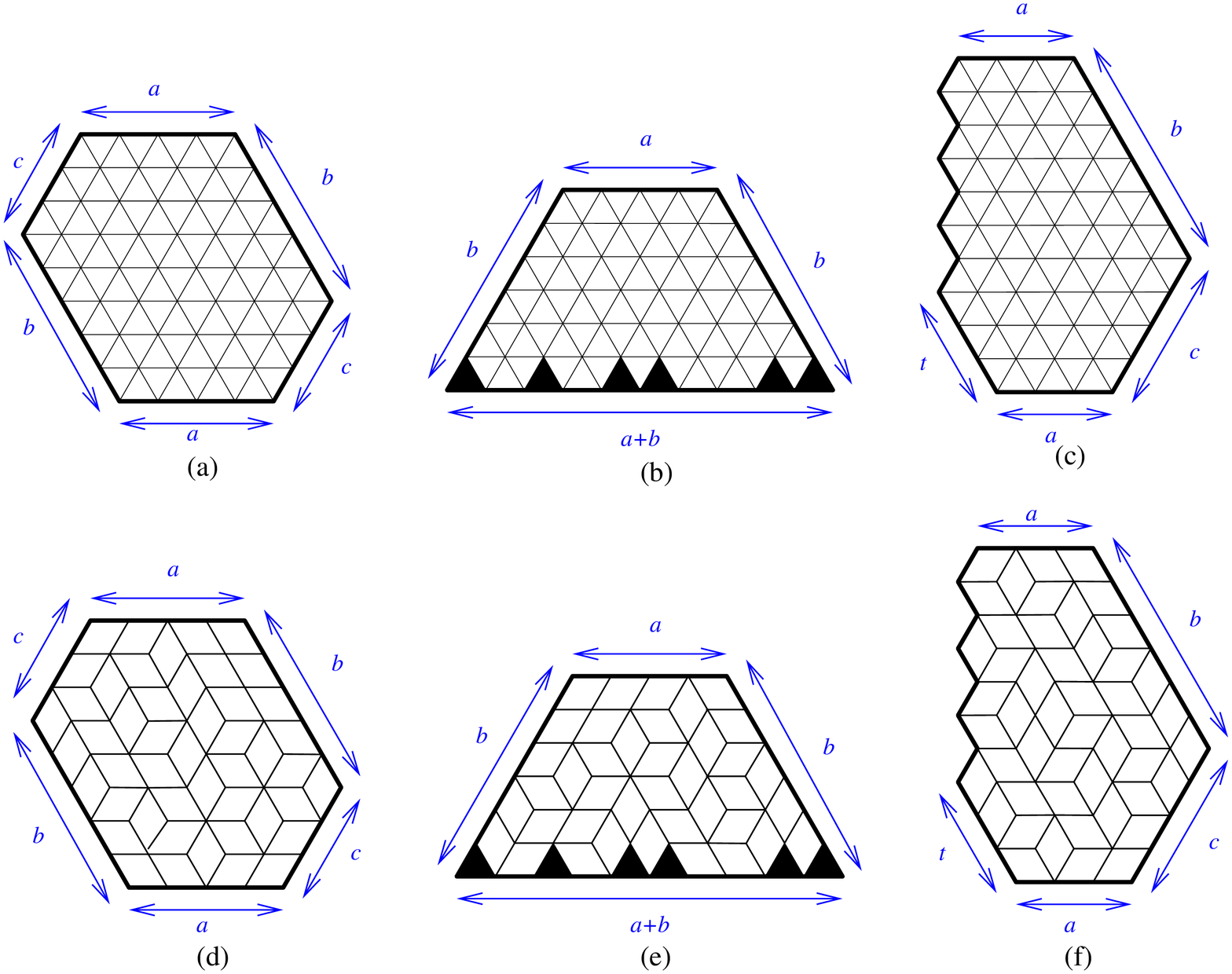}
\caption{Three popular regions in the enumeration of lozenge tilings: (a) the quasi-regular hexagon, (b) the semi-hexagon with dents on the base, and (c) the halved hexagon. Sample tilings are shown in the second row.}\label{Fig:3polar}
\end{figure}

We will frequently mention the following three well-known tiling enumerations of the centrally symmetric hexagon (see Figure \ref{Fig:3polar}(a)), the semi-hexagon with dents on the base (see Figure \ref{Fig:3polar}(b)), and the halved hexagon (illustrated in Figure \ref{Fig:3polar}(c)). The enumeration of tilings of the hexagon of side-lengths $a,b,c,a,b,c$ (in counter-clockwise order, starting from the north side\footnote{From now on, we always list the side-lengths of hexagonal regions in this order.}) is credited to P. A. MacMahon \cite{Mac} in the early 1900s. However, MacMahon did not work on tilings; what he proved is a more general result (see Theorem \ref{MacMahon2}) on the enumeration of plane partitions fitting in a given box (or `\emph{boxed plane partitions}').  The tiling enumeration of the semi-hexagon with dents on the base is due to H. Cohn, M. Larsen, and J. Propp \cite[Proposition 2.1]{CLP} when they give a bijection between lozenge tilings of the region and the semi-strict Gelfand--Tsetlin patterns \cite{GT}. It is well-known that the lozenge tilings of a dented semi-hexagon also correspond to the \emph{column-strict plane partitions} (or \emph{reverse semi-standard Young tableaux}). The enumeration of the halved hexagon was first found by R. Proctor \cite[Corollary 4.1]{Proc} in the form of the number of a particular class of staircase plane partitions. His result also implies the enumeration of the \emph{transpose-complementary plane partitions}, one of the ten symmetry classes of plane partitions \cite{Stanley2}.

\subsection{Generalized boxed plane partitions}
\begin{figure}\centering
\includegraphics[width=15cm]{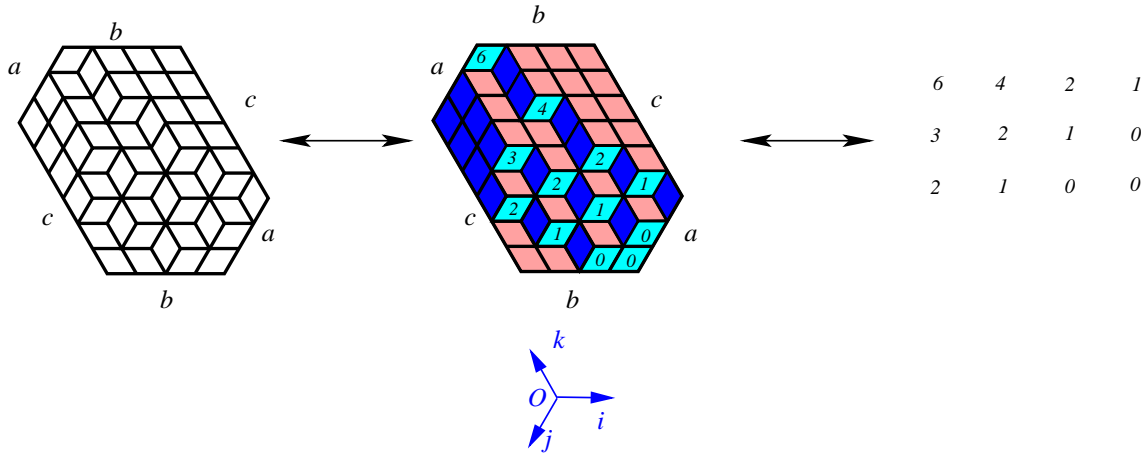}
\caption{Correspondence between lozenge tilings of a hexagon, stacks of unit cubes fitting in a rectangular box, and plane partitions. The picture first appeared in \cite{Tri18}.}\label{Macfig}
\end{figure}

A \emph{plane partition} can be defined as a rectangular array of non-negative integers with weakly decreasing rows and columns.  One can view a plane partition as a monotonic stack of unit cubes fitting in a given rectangular box, and the latter, in turn, are in bijection with lozenge tilings of a quasi-regular hexagon \cite{DT}. For example, we can write the entries of the plane partition $\pi$ in the right picture of Figure \ref{Macfig} on a rectangular board of the same size (in this case, a $3\times 4$ board) embedded on the plane $Oij$, and we place the corresponding number of unit cubes on each entry of the board. This way, one can interpret the plane partition $\pi$ as a monotonic stack of unit cubes in the middle picture of Figure \ref{Macfig}. This stack, in turn, can be projected on the plane $i+j+k=0$ to obtain the lozenge tiling of a hexagon shown in the left picture. From this point of view, MacMahon's classical theorem \cite{Mac} on boxed plane partitions can be stated in the language of the volume generating function of the stacks as follows.

 Let $q$ be an indeterminate. The \emph{$q$-factorial} is defined as $[n]_q!:=\prod_{i=1}^{n}\frac{1-q^n}{1-q}$, where $[0]_q!=1$; and the \emph{$q$-hyperfactorial} is $\Hf_q(n):=[0]_q![1]_q!\dotsc[n-1]_q!$, where $\Hf_q(0)=1$. 
 \begin{thm}[MacMahon's Theorem \cite{Mac}]\label{MacMahon2} For non-negative integers $a,b,c$
\begin{align}\label{Maceq}
\sum_{\pi}q^{|\pi|}&=\prod_{i=1}^{a}\prod_{j=1}^{b}\prod_{k=1}^{c}\frac{q^{i+j+k-1}-1}{q^{i+j+k-2}-1}\notag\\
&=\frac{\Hf_q(a)\Hf_q(b)\Hf_q(c)\Hf_q(a+b+c)}{\Hf_q(a+b)\Hf_q(b+c)\Hf_q(c+a)},
\end{align}
where the sum is taken over all monotonic stacks $\pi$ fitting in an $(a\times b \times c)$-box, and where $|\pi|$ denotes the \emph{volume}  of $\pi$.
\end{thm}

Setting $q$ tend to 1, one obtains the tiling number of the quasi-regular hexagon from the above theorem. 

The beauty of formula (\ref{Maceq})  has inspired a large body of work, focusing on the enumeration of lozenge tilings of hexagons with defects. Put differently, MacMahon's theorem gives a weighted enumeration of lozenge tilings of a hexagon. As the enumeration is a function in $q$, we often call it a ``\emph{$q$-enumeration}." Unfortunately, such elegant $q$-enumerations are \emph{very rare} in the domain of enumeration of lozenge tilings. Besides the three new $q$-enumerations below, only a few are known (see, e.g., \cite{Stanley,Stanley2,Krat3} and the list of references therein).

\begin{figure}\centering
\includegraphics[width=14cm]{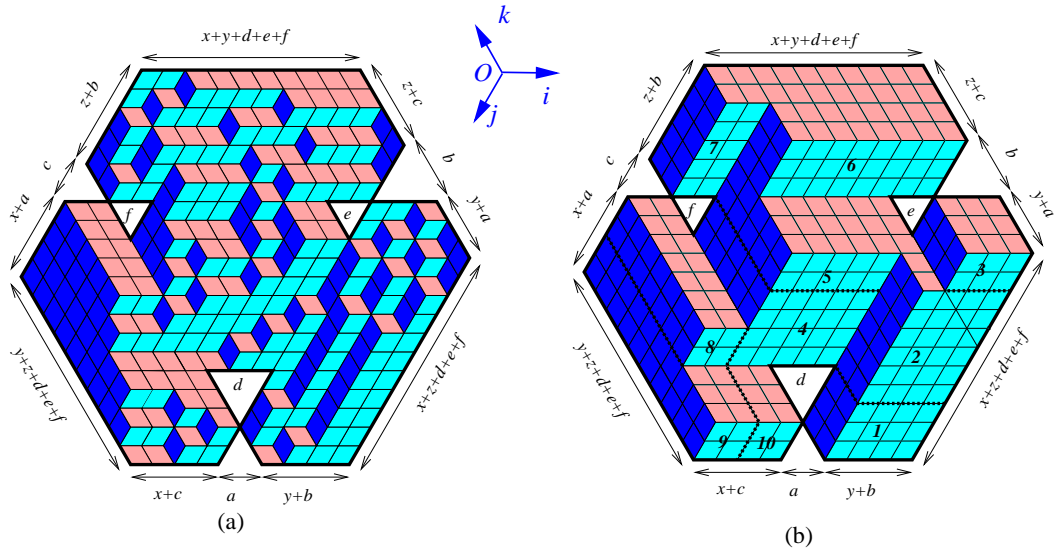}
\caption{(a) Viewing a lozenge tiling of a hexagon with three dents as a stack of unit cubes fitting in a compound box. (b) The empty compound box with the floors of the rooms labeled by $1,2,\dotsc,10$. The picture first appeared in \cite{Tri19}.}
\label{NSFbox2}
\end{figure}

In \cite{Tri18,Tri19,Twofern}, MacMahon's classical Theorem \ref{MacMahon2} has been generalized by investigating $q$-enumerations of several families of hexagons with defects on the boundary. We usually call this type of defects ``\emph{dents}." In order to give an illustrated example, let us focus on only the main result in \cite{Tri19}.

Our region here is a $9$-parameter generalization of the dented hexagon that first appeared in Problem 3\footnote{This problem was first solved by T. Eisenk\"{o}lbl \cite{Eisen}.} on Propp's well-known list of open problems \cite{Propp}. In particular, we consider a certain hexagon with three bowtie-shaped dents on three non-consecutive sides, denoted by $F\begin{pmatrix}x&y&z\\a&b&c\\d&e&f\end{pmatrix}$. We view the lozenge tilings of a $F$-type region as monotonic stacks of unit cubes fitting in a \emph{compound box} $\mathcal{B}=\mathcal{B}\begin{pmatrix}x&y&z\\a&b&c\\d&e&f\end{pmatrix}$, which is the union of 10 `rooms' (see Figure \ref{NSFbox2}(a)).  In Figure \ref{NSFbox2}(b), we have a picture of the lozenge tiling corresponding to the empty stack, and this also gives a $3$-D picture of the compound box $\mathcal{B}$. The floors of the rooms are labelled by $1,2,\dotsc,10$. One readily sees that the stacks of unit cubes here have the same monotonicity as the ordinary plane partitions (namely, the tops of the columns are weakly decreasing along the directions of the $i$- and $j$-axes). Like the case of MacMahon's theorem, the volume generating function of these stacks is always given by a simple product formula in terms of $q$-hyperfactorials.

\begin{thm}[Theorem 1.2 in \cite{Tri19}]\label{Bowtiethm} For non-negative integers $a,$ $b,$ $c,$ $d,$ $e,$ $f,$ $x,$ $y,$ $z$
\small{\begin{align}\label{maineqbox}
& \sum_{\pi}q^{|\pi|}=\notag\\
 &\frac{\Hf_q(x)\Hf_q(y)\Hf_q(z)\Hf_q(a)^2\Hf_q(b)^2\Hf_q(c)^2\Hf_q(d)\Hf_q(e)\Hf_q(f)\Hf_q(d+e+f+B)^4}{\Hf_q(a+d)\Hf_q(b+e)\Hf_q(c+f)\Hf_q(d+e+B)\Hf_q(e+f+B)\Hf_q(f+d+B)}\notag\\
 &\times \frac{\Hf_q(A+2B)\Hf_q(A+B)^2}{\Hf_q(A+B+x)\Hf_q(A+B+y)\Hf_q(A+B+z)}\notag\\
  &\times \frac{\Hf_q(a+b+d+e+B)\Hf_q(a+c+d+f+B)\Hf_q(b+c+e+f+B)}{\Hf_q(a+d+e+f+B)^2\Hf_q(b+d+e+f+B)^2\Hf_q(c+d+e+f+B)^2}\notag\\
  &\times\frac{\Hf_q(a+d+x+y)\Hf_q(b+e+y+z)\Hf_q(c+f+z+x)}{\Hf_q(a+b+y)\Hf_q(b+c+z)\Hf_q(c+a+x)}\notag\\
  &\times \frac{\Hf_q(A-a+B+z)\Hf_q(A-b+B+x)\Hf_q(A-c+B+y)}{\Hf_q(b+c+e+f+B+z)\Hf_q(c+a+d+f+B+x)\Hf_q(a+b+d+e+B+y)},
\end{align}}
\normalsize where the sum is taken over all monotonic stacks $\pi$ fitting in the  compound box $\mathcal{B}=\mathcal{B}\begin{pmatrix}x&y&z\\a&b&c\\d&e&f\end{pmatrix}$, and where $A=a+b+c+d+e+f$, $B=x+y+z$.
\end{thm}

We believe that there are more elegant $q$-enumerations of lozenge tilings waiting for us to explore. 
\begin{prob}\label{problem1}
Find more dented regions whose corresponding volume generating functions are given by simple product formulas.
\end{prob}
We hope that, after collecting enough examples of nice $q$-enumerations, we can solve the following problem:
\begin{prob}\label{problem2}
Characterize the compound boxes, which yield nice volume generating functions.
\end{prob}

It is worth noticing that the above ``stack-box model" does \emph{not} work well when our regions have ``\emph{holes}" (i.e., some portions removed from its interior). In this case, there may be more than one way to lift a lozenge tiling to a stack of unit cubes.

\subsection{Elliptic Weight}

A. Borodin, V. Gorin, and E. M. Rains \cite{Borodin} provide a different way to define weight for tilings, called ``elliptic weight," as follows. Then the weight of a tiling is the product of the weights of its lozenges. We will adapt and specialize Borodin--Gorin--Rains' elliptic weight in this section.

\begin{figure}\centering
\includegraphics[width=8cm]{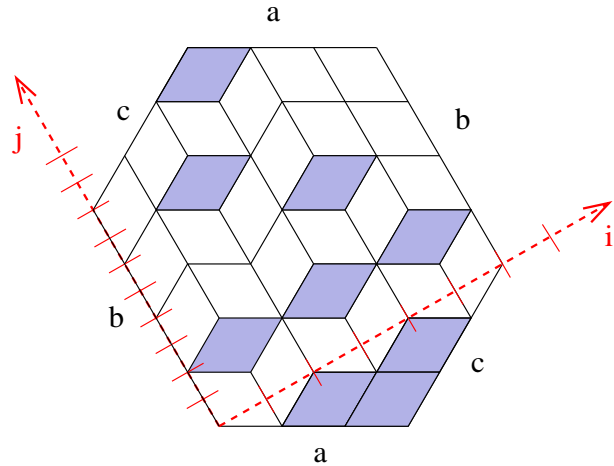}
\caption{Assigning weight to lozenges.}\label{Fig:weight}
\end{figure}

The $j$-axis of our coordinate runs along  a lattice line with the unit equal to $1/2$ times the side-length of a lozenge. The $i$-axis is perpendicular with the $j$-axis at a lattice vertex (this vertex is the origin of our coordinate system); the unit on the $i$-axis is equal to $\sqrt{3}/2$ times the side-length of a lozenge. Figure \ref{Fig:weight} shows a particular placement of our coordinate system.

Only one of the three types of lozenges, as shown in Figure \ref{Fig:orientation}, have the diagonals parallel to the $i$- and $j$-axes. (In Figure \ref{Fig:weight}, these lozenges are the right lozenges.)
Each lozenge of this type with center at the point $(i,j)$ is weighted by
\begin{equation}
\wt_1(i,j)=\frac{Xq^{j}+Yq^{-j}}{2},
\end{equation}
where $X,Y,q$ are three indeterminates.  (The weight does not depend on $i$.) All other lozenges have weight $1$. The weight of a tiling is now the product of weights its lozenges.  In the rest of the paper, we use the notation $\M(R)$ for the sum of weights of all tilings of a weighted region $R$. We call $\M(R)$ the \emph{tiling generating function} of $R$. When $R$ is unweighted, $\M(R)$ is exactly the tiling number of $R$.

This weight behaves very well and could be considered as the generalization of the ``volume weight" investigated in the previous section. Indeed, when $X=2$ and $Y=0$, the weight becomes $\wt_2(i,j)=q^{j}$. With the weighting system as in Figure \ref{Fig:weight}, each tiling $\tau$ of the hexagon is weighted by $C\cdot q^{2\cdot Vol(\tau)}$, where $Vol(\tau)$ is the volume of the stack of unit cubes corresponding to the tiling $\tau$, and $C$ is a constant independent from the choice of the tiling $\tau$. Sometimes, we call $\wt_2$ the  ``natural weight" (as it is essentially equivalent to the weight of the boxed plane partitions).

We are also interested in the following specialization of $\wt_1$:
\begin{equation}
\wt_3(i,j)=\frac{q^{j}+q^{-j}}{2}.
\end{equation}
One can view $\wt_3$ as a symmetrization of $\wt_2$, and we often call $\wt_3$ the ``symmetric weight." As mentioned in the previous section, the natural weight $\wt_2(i,j)=q^{j}$ does not often give nice tiling generating functions. In contrast, the symmetric weight $\wt_3(i,j)$ behaves much better. For example, it has been shown by M. Ciucu, T. Eisenk\"{o}lbl, C. Krattenthaler, and D. Zare \cite{CEKZ} that the ``plain" tiling number (unweighted counting of tilings) of a ``\emph{cored hexagon}" is always given by a simple product formula (see Figure \ref{Fig:weightn}). However, there is no such formula for the tiling generating function associated with $\wt_2$. On the other hand, as shown by H. Rosengren \cite{Rosen}, the symmetric weight $\wt_3$ gives a simple product formula for the tiling generating function of the cored hexagon. We have observed a similar fact for a halved hexagon with defects \cite{LR20}: the natural weight $\wt_2$ does not give a nice $q$-enumeration, but the symmetric weight $\wt_3$ does. 

Despite its very nice behavior, the study about the weight $\wt_3$ is still \emph{extremely} limited.  It deserves more attention from experts in the field. It would be interesting to see if the weight $\wt_3$ yields nice tiling generating functions for known regions.

\begin{prob}\label{problem3}
Find the tiling generating functions with respect to the symmetric  weight $\wt_3$ for known families of regions. For instance, one would like to find the tiling generating function for the ``\emph{$S$-cored hexagon}" introduced by Ciucu and Krattenthaler \cite{CK13} as a generalization of the cored hexagon in \cite{CEKZ} (see Figure \ref{Fig:offcenter} for an example of S-cored hexagon). This would give a generalization for Rosengren's enumeration in \cite{Rosen}.
\end{prob}

As most of the known tiling enumerations are unweighted ones, there would be many things to do with this problem.

\medskip

All the weights $\wt_1,\wt_2,\wt_3$ can be viewed as some rational functions in $z=q^{j}$. We want to find more weights of this type, such that they give nice tiling  generating functions, say for several well-known families of regions, like the hexagons, the semi-hexagons, and  the halved hexagons.

\begin{prob}\label{problem4}
Find new ``nice'' weights that are rational functions in $q^{j}$.
\end{prob}

\begin{figure}\centering
\begin{picture}(0,0)%
\includegraphics{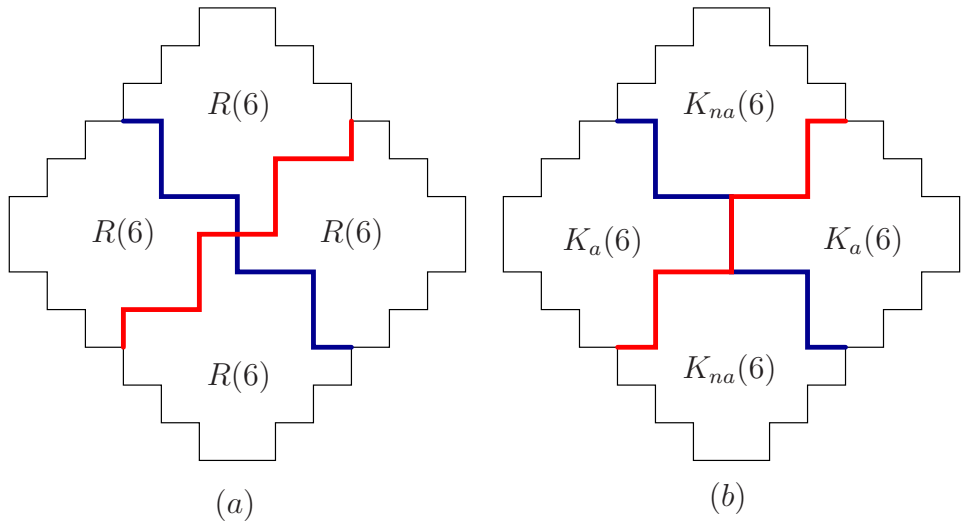}%
\end{picture}%
%
%
\setlength{\unitlength}{3947sp}%
\begingroup\makeatletter\ifx\SetFigFont\undefined%
\gdef\SetFigFont#1#2#3#4#5{%
  \reset@font\fontsize{#1}{#2pt}%
  \fontfamily{#3}\fontseries{#4}\fontshape{#5}%
  \selectfont}%
\fi\endgroup%
\begin{picture}(5929,3252)(1170,-2873)
\put(2458,-2798){\makebox(0,0)[lb]{\smash{{\SetFigFont{12}{14.4}{\rmdefault}{\mddefault}{\updefault}{$(a)$}%
}}}}
\put(5529,-2751){\makebox(0,0)[lb]{\smash{{\SetFigFont{12}{14.4}{\rmdefault}{\mddefault}{\updefault}{$(b)$}%
}}}}
\put(2410,-295){\makebox(0,0)[lb]{\smash{{\SetFigFont{12}{14.4}{\rmdefault}{\mddefault}{\updefault}{$R(6)$}%
}}}}
\put(3109,-1104){\makebox(0,0)[lb]{\smash{{\SetFigFont{12}{14.4}{\rmdefault}{\mddefault}{\updefault}{$R(6)$}%
}}}}
\put(2410,-1995){\makebox(0,0)[lb]{\smash{{\SetFigFont{12}{14.4}{\rmdefault}{\mddefault}{\updefault}{$R(6)$}%
}}}}
\put(1692,-1104){\makebox(0,0)[lb]{\smash{{\SetFigFont{12}{14.4}{\rmdefault}{\mddefault}{\updefault}{$R(6)$}%
}}}}
\put(5366,-295){\makebox(0,0)[lb]{\smash{{\SetFigFont{12}{14.4}{\rmdefault}{\mddefault}{\updefault}{$K_{na}(6)$}%
}}}}
\put(5366,-1954){\makebox(0,0)[lb]{\smash{{\SetFigFont{12}{14.4}{\rmdefault}{\mddefault}{\updefault}{$K_{na}(6)$}%
}}}}
\put(6238,-1145){\makebox(0,0)[lb]{\smash{{\SetFigFont{12}{14.4}{\rmdefault}{\mddefault}{\updefault}{$K_a(6)$}%
}}}}
\put(4620,-1151){\makebox(0,0)[lb]{\smash{{\SetFigFont{12}{14.4}{\rmdefault}{\mddefault}{\updefault}{$K_a(6)$}%
}}}}
\end{picture}%

\caption{Three kinds of quartered Aztec diamonds of order 6. The figure first appeared in \cite{Tri6}.}
\label{Fig.QA}
\end{figure}

\begin{figure}\centering
\includegraphics[width=10cm]{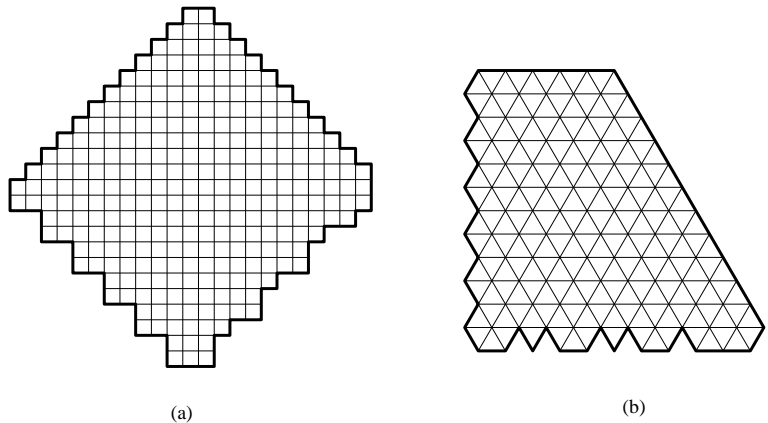}
\caption{(a) A quartered Aztec rectangle and (b) A quartered hexagon.}\label{QHsurvey}
\end{figure}

Jockusch and Propp \cite{JP} introduced the ``\emph{quartered Aztec diamonds}" as quarters of an Aztec diamond divided by two zigzag cuts passing the center (see Figure \ref{Fig.QA}). These regions have been re-investigated and generalized in \cite{Trinewquarter,Trinewquarter2,Trinewquarter3,Tri6}. These papers showed that one could transform a ``\emph{quartered Aztec rectangle}" (a natural generalization of the quartered Aztec diamond) into a quartered hexagons using certain local graph transformations. See Figure \ref{QHsurvey} for an example of a quartered Aztec rectangle and a quartered hexagon. It turns out the tiling numbers of the two regions are only different by a multiplicative factor, which is a perfect power of $2$. As a nice $q$-formula for tiling generating function of  the quartered hexagon has been found in \cite{LR20}, one would like to find such a $q$-formula for the quartered Aztec diamond and quartered Aztec rectangles.

\begin{prob}\label{problem5}
Find a nice $q$-enumeration for the domino tilings of the quartered Aztec diamond and quartered Aztec rectangle.
\end{prob}

\subsection{Unusual weights}

A special weight of lozenge tilings of a hexagon has been inspired by the well-known trace formula of R. Stanley. By letting $b\rightarrow \infty$ in MacMahon's classical formula in  Theorem \ref{MacMahon2}, we obtain another well-known formula of MacMahon:
\begin{equation}\label{Macformula2}
\sum_{\pi\in \mathcal{P}(a,c)}q^{|\pi|}=\prod_{i=1}^{a}\prod_{i=1}^{c}\frac{1}{1-q^{i+j-1}},
\end{equation}
where $\mathcal{P}(a,c)$ denotes the set of plane partitions with at most $c$ rows and $a$ columns. Generalizing MacMahon's formula (\ref{Macformula2}), Stanley \cite{Stanley3}  proved the \emph{trace generating function}:
\begin{equation}
\sum_{\pi\in\mathcal{P}(a,c)}q^{|\pi|}t^{\tr(\pi)}=\prod_{i=1}^{a}\prod_{i=1}^{c}\frac{1}{1-tq^{i+j-1}},
\end{equation}
where the \emph{trace} is defined by $\tr(\pi):=\sum_{i} \pi_{i,i}$.
E. R. Gansner  \cite{Gan1,Gan2} later extended Stanley's work by showing that
\begin{equation}
\sum_{\pi\in\mathcal{P}(a,c)}\prod_{-c<\ell<a}q_{\ell}^{\tr_{\ell}(\pi)}=\prod_{i=0}^{a-1}\prod_{j=0}^{c-1}\left(1-\prod_{\ell=-i}^{j}q_{\ell}\right)^{-1},
\end{equation}
where the \emph{$\ell$-trace} is defined as $\tr_{\ell}(\pi):=\sum_{j-i=\ell}\pi_{i,j}$.

Strictly speaking, Stanley's and Gansner's trace formulas above do not give any weighted enumerations for boxed plane partitions (equivalently, lozenge tilings of hexagons), as there is no upper bound for the parts of plane partitions. Recently, S. Kamioka   \cite{Kam} provides elegant boxed versions for the above formulas.

\begin{figure}\centering
\setlength{\unitlength}{3947sp}%
\begingroup\makeatletter\ifx\SetFigFont\undefined%
\gdef\SetFigFont#1#2#3#4#5{%
  \reset@font\fontsize{#1}{#2pt}%
  \fontfamily{#3}\fontseries{#4}\fontshape{#5}%
  \selectfont}%
\fi\endgroup%
\resizebox{!}{4.6cm}{
\begin{picture}(0,0)%
\includegraphics{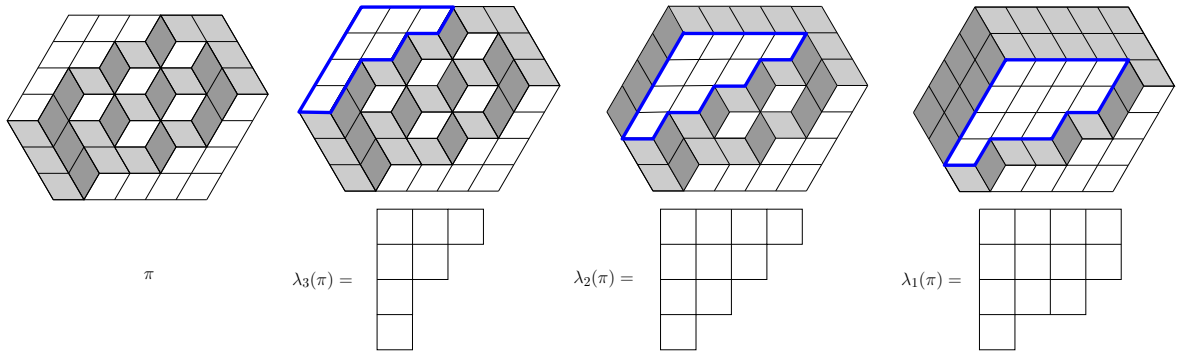}%
\end{picture}%
%
%

\begin{picture}(15588,4663)(1200,-6260)
\put(13131,-5371){\makebox(0,0)[lb]{\smash{{\SetFigFont{15}{14.4}{\rmdefault}{\mddefault}{\updefault}{\color[rgb]{0,0,0}$\lambda_1(\pi)=$}%
}}}}
\put(5011,-5371){\makebox(0,0)[lb]{\smash{{\SetFigFont{15}{14.4}{\rmdefault}{\mddefault}{\updefault}{\color[rgb]{0,0,0}$\lambda_3(\pi)=$}%
}}}}
\put(8761,-5351){\makebox(0,0)[lb]{\smash{{\SetFigFont{15}{14.4}{\rmdefault}{\mddefault}{\updefault}{\color[rgb]{0,0,0}$\lambda_2(\pi)=$}%
}}}}
\put(3036,-5302){\makebox(0,0)[lb]{\smash{{\SetFigFont{15}{14.4}{\rmdefault}{\mddefault}{\updefault}{\color[rgb]{0,0,0}$\pi$}%
}}}}
\end{picture}}
\caption{The plane partition $\pi$ in the form of a stack of unit cubes and three cross-sections of $\pi$.}\label{Fig:crosssection}
\end{figure}

For a plane partition $\pi$ of shape $\lambda$ and a number $1\leq k\leq \pi_{1,1}$, we define the  \emph{$k$-truncation} $\pi^{(k)}$ of $\pi$ to be the plane partition obtained from $\pi$ by removing all entries less than $k$. The shape $\lambda^{(k)}(\pi)$ of $\pi^{(k)}$ is called the \emph{$k$-cross-section} of $\pi$. In particular, we have $\pi^{(1)}=\pi$ and $\lambda^{(1)}(\pi)=\lambda$ (see Figure \ref{Fig:crosssection} for the interpretation of the cross-sections). The \emph{Durfee square} of a partition  is the largest square fitting in its Ferrers diagram.

We often use the standard  $q$-Pochhammer symbol in our tiling formulas:
\begin{equation}
(x; q)_n:=
\begin{cases}1 &\text{if $n=0$;}\\
(1-x)(1-xq)\cdots(1-xq^{n-1}) &\text{if $n>0$;}\\
\frac{1}{(1-xq^{-1})(1-xq^{-2})\cdots(1-xq^{n})} &\text{if $n<0$.}\\
\end{cases}
\end{equation}
Strictly speaking, the   above $q$-Pochhammer  symbol is not well-defined when $n$ is a negative integer and $a=q^{k}$ for some $1\leq k \leq -n$. However, this is not the case in our paper.

\begin{thm} [Theorem 9 in \cite{Kam}] \label{Kamiokathm1}
\begin{equation}\label{Kamiokaeq}
\sum_{\pi\in \mathcal{P}(a,b,c)}q^{|\pi|}t^{\tr(\pi)}\prod_{k=1}^{\pi_{1,1}}\frac{(q^{n-k+1};q)_{D_{k}(\pi)}}{(tq^{n-k+1};q)_{D_{k}(\pi)}}
=\prod_{i=1}^{a}\prod_{j=1}^{b}\prod_{k=1}^{c}\frac{1-tq^{i+j+k-1}}{1-tq^{i+j+k-2}},\end{equation}
where $D_{k}(\pi)$ is the side-length of the Durfee square of the $k$-cross-section  $\lambda^{(k)}(\pi)$ of $\pi$.
\end{thm}
This enumeration can be viewed as a weighted enumeration of the lozenge tilings of a hexagon. In particular, the tiling $\tau=\tau_{\pi}$ corresponding to the plane partition $\pi$ is weighted by $\wt_K(\tau)=q^{|\pi|}t^{\tr(\pi)}\prod_{k=1}^{\pi_{1,1}}\frac{(q^{n-k+1};q)_{D_{k}(\pi)}}{(tq^{n-k+1};q)_{D_{k}(\pi)}}$. It is worth noticing that Kamioka also proved a finite version of Gansner's formula (see Theorem 17 in the same paper). He later generalizes further the results to arbitrary shapes in \cite{Kam2}.

\begin{figure}\centering
\includegraphics[width=7cm]{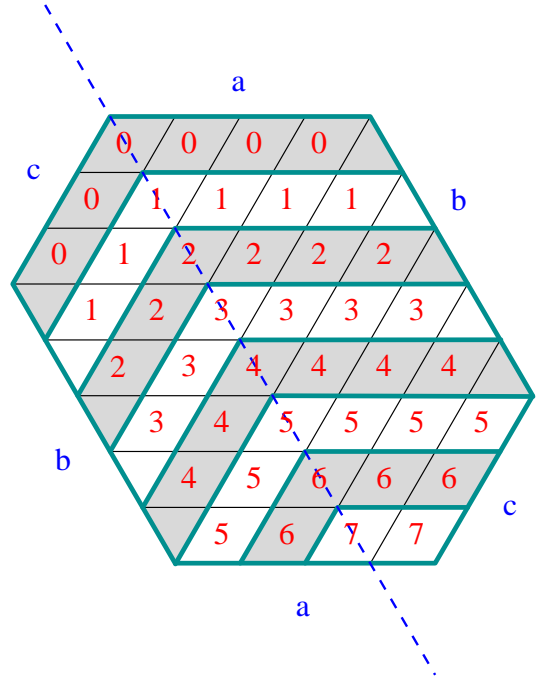}
\caption{The weight for lozenge tiling inspired by Kamioka's formula.}\label{Fig:Kamiokaweight}
\end{figure}

Unlike the weights $\wt_1,\wt_2,\wt_3$ in the previous section, Kamioka's theorem does \emph{not} provide the weights for lozenges. However, one could define  the corresponding  lozenge-weights  as follows.
Divide the quasi-regular hexagon $Hex(a,b,c)$ into hooks as in Figure \ref{Fig:Kamiokaweight}. Each right lozenge is labeled as in the figure.  
Now, each right lozenge with label $x$ is weighted by $q^{-x}$, except for the ones intersecting with the dotted line, which are weighted by $(tq)^{-x}\frac{(q^{n};q^{-1})_{n-x}}{(tq^{n};q^{-1})_{n-x}}$. All lozenges of different orientations are weighted by $1$. It is easy to see that the product of lozenge-weights in tiling $\tau_{\pi}$ is equal to $C\cdot \wt_K(\tau_{\pi})$, where $C$ is a constant independent from  the choice of tilings. We get back Kamioka's weight by normalizing the lozenge-weights. See  \cite{Trace} for more details.

In his classical paper \cite{Stanley2}, Stanley lists ten symmetry classes of plane partitions. Each of the ten classes is equivalent to a particular type of symmetric tilings of a hexagon. We would expect the existence of Kamioka's trace formula for symmetric tilings.

\begin{prob}\label{problem6}
Find versions of Kamioka's trace formula for symmetric tilings.
\end{prob}

\begin{figure}\centering
\includegraphics[width=13cm]{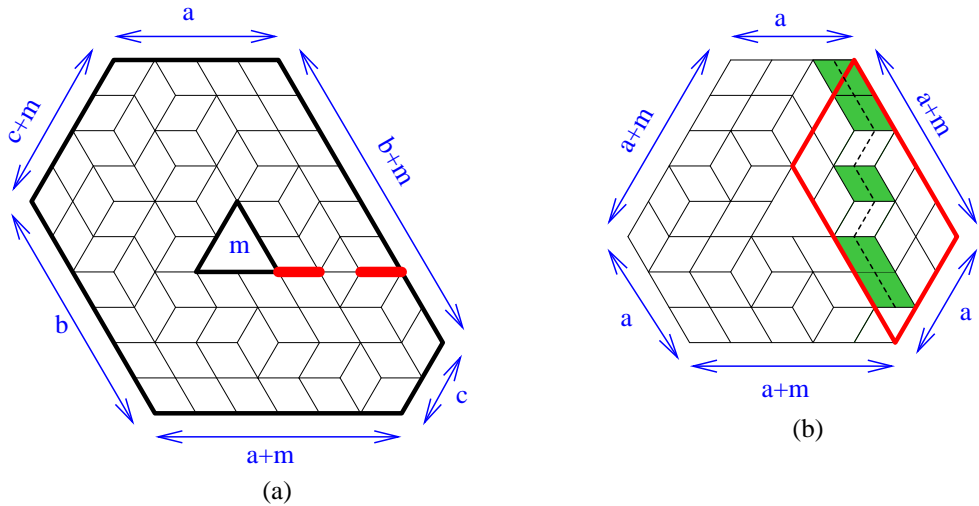}
\caption{Weight in the (-1)-enumerations of a cored hexagon.}\label{Fig:weightn}
\end{figure}

We conclude this subsection by investigating another unusual weight of tilings inspired by the \emph{(-1)-phenomenon} found by J. Stembridge \cite{Stem2}, and more generally, the \emph{Cyclic Sieving Phenomenon} by V. Reiner, D. Stanton, and D. White \cite{RSW}. In \cite{CEKZ}, the authors proved striking formulas for specific $(-1)$-enumerations of the tilings of  the cored hexagons and its cyclically symmetric tilings (see Theorems 4, 5, and 7 therein). In particular, if we extend the base of the triangular hole in the cored hexagon to the right, then each tiling $\tau$ is weighted by $(-1)^{n(\tau)}$, where $n(\tau)$ is the number of edges of lozenges of the tiling $\tau$ contained in the extended side (see Figure \ref{Fig:weightn}(a); in this case $n(\tau)=2$). The critical point of the definition of $n(\tau)$ is that when we encode each tiling as a family of the non-intersecting lattice paths in the spirit of the well-known Lindstr\"{o}m--Gessel--Viennot Theorem (see, e.g., \cite{Lind, Stem3, GV}), the sign of the path family is exactly $(-1)^{n(\tau)}$. For the case of cyclically symmetric tilings, the weight is quite different. We consider the fundamental region of the cyclically symmetric tilings (which is illustrated as the parallelogram with bold sides in Figure \ref{Fig:weightn}(b)), then each tiling is weighted by $(-1)^{n_6(\tau)}$, where $n_6(\tau)$ is the sum of the horizontal distances between the shaded lozenges and the lower-left border of the fundamental region ($n_6(\tau)=0+0+1+2+2=5$ in Figure \ref{Fig:weightn}).

One would ask for similar $(-1)$-enumerations for the following three generalizations of the cored hexagons: 
\begin{enumerate}
\item[(1)] the ``\emph{$S$-cored hexagons}" (i.e., a hexagon with a cluster of four triangles removed from the center; see Figure \ref{Fig:Shamrock}(a)) in \cite{CK13}, 
\item[(2)] the ``\emph{$F$-cored hexagons}" (i.e., a hexagon with an array of alternating triangles removed from the center) in \cite{Ciu1}, and
\item[(3)] the hexagons with three arrays of triangles removed in \cite{HoleDent}. 
\end{enumerate}
It is worth noticing that the unweighted enumeration of cyclically symmetric tilings of a $S$-cored hexagon has been provided by Ciucu in \cite{cyclicshamrock}.

\begin{prob}\label{problem7}
(a) Find the $(-1)$-enumerations of tilings of the three generalizations of core hexagons listed above.

(b) Find the $(-1)$-enumeration of cyclically symmetric tilings for the $S$-cored hexagon.
\end{prob}

Recall that Rosengren proved a simple product formula for the tiling generating of the cored hexagons using the weight $\wt_3=\frac{q^{j}+q^{-j}}{2}$. We note that the versions of Rosengren's $q$-enumeration for the three generalizations of the cored hexagons are not known at this point.  

\begin{prob}\label{problem8}
 Find a version of Rosengren's $q$-enumeration for each of the three generalizations of the cored hexagons above.
 \end{prob}

Inspired by the above $(-1)$-enumeration, we would like to find a signed version of Rosengren's $q$-enumeration for the cored hexagons, and more generally, for its  generalizations.

\begin{prob}\label{problem9}
(a) Find a signed version of Rosengren's $q$-enumeration for the cored hexagon.

(b) Generalize part (a) to the three generalizations of the cored hexagons above.
\end{prob}

\section{Shuffling Phenomenon}

This section is devoted to an exciting property of tilings, the  ``shuffling phenomenon"\footnote{The ``shuffling phenomenon" here is \emph{not} related to the ``domino shuffling" operation in \cite{Elkies1,Elkies2}.}, that was first introduced in \cite{shuffling}. The minor modification in a region would lead to an unpredictable change of the tiling number. However, in some particular situations, the tiling number is changed by only  a simple multiplicative factor. See e.g. \cite{shuffling2, shuffling3, Byun, Byun2, Ful, Ful2, CLR, Con, Semitwodent, LR20} for recent development of the phenomenon.

Let $x,y,z,u,d$ be nonnegative integers, such that $u,d \leq n$. Consider a symmetric hexagon of side-lengths $x+n-u,y+u,y+d,x+n-d,y+d,y+u$. We remove $u$ up-pointing and $d$ down-pointing unit triangles along the  horizontal lattice line $l$ that contains the west and the east vertices of the hexagon. Let $U=\{s_1,s_2,\dotsc,s_u\}$ and $D=\{t_1,t_2,\dots, t_d\}$ denote the position sets of the up-pointing and down-pointing removed unit triangles (ordered from left to right), respectively. Assume that $n$ is the size of the union $U\cup D$.  Denote by $H_{x,y}(U;D)$ the hexagon with the above setup of removed unit triangles. We call it a \emph{doubly--dented hexagon}. We now allow to `shuffle' the positions of the up- and down-pointing unit triangles in the symmetric difference $U\Delta D$ to obtain new position sets $U'$ and $D'$, respectively. The following theorem shows that the above shuffling unit triangles changes the tiling number of the region by only a simple multiplicative factor. Moreover, the factor can be written in a similar form to Cohn--Larsen--Propp's  tiling formula of a semi-hexagon \cite[Proposition 2.1]{CLP}.

\begin{figure}\centering
\includegraphics[width=13cm]{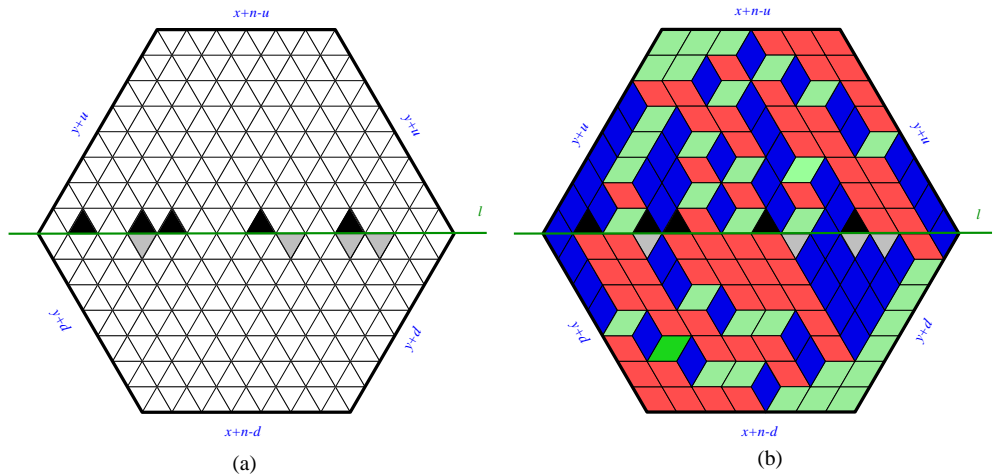}
\caption{The region $H_{4,3}(2,4,5,8,11;\ 4,9,11,12)$ (left) and a lozenge tiling of its (right). The back and shaded triangles indicate the unit triangles removed. }\label{multiplefernfig}
\end{figure}

\begin{thm}[Shuffling Theorem; Theorem 2.1 in \cite{shuffling}]\label{factorthm} For nonnegative integers $x,y,u,d,n$ ($u,d \leq n$) and four ordered subsets $U=\{s_1,s_2,\dotsc,s_u\}$,  $D=\{t_1,t_2,\dots, t_d\}$, $U'=\{s'_1,s'_2,\dotsc,s'_{u}\}$, and  $D'=\{t'_1,t'_2,\dots, t'_{d}\}$ of $[x+y+n]$,  such that $U\cup D =U'\cup D'$ and $U\cap D =U'\cap D'$. 
 We have
\begin{equation}\label{genmaineq}
  \frac{\M(H_{x,y}(U;D))}{\M(H_{x,y}(U';D'))}= \frac{\displaystyle\prod_{1\leq i <j\leq u}\frac{s_j-s_i}{j-i}\cdot\displaystyle \prod_{1\leq i <j\leq d}\frac{t_j-t_i}{j-i}}{\displaystyle\prod_{1\leq i <j\leq u}\frac{s'_j-s'_i}{j-i}\cdot\displaystyle\prod_{1\leq i <j\leq d}\frac{t'_j-t'_i}{j-i}},
\end{equation}
recall that we use the notation $\M(R)$ for the number of tilings of the unweighted region $R$ (when $R$ is weighted, $\M(R)$ denotes the tiling generating function of $R$).
\end{thm}
We would like to emphasize that neither $H_{x,y}(U;D)$ nor $H_{x,y}(U;'D')$ has a nice tiling number in general. We also note that several stronger versions of the above Shuffling Theorem were provided in \cite{shuffling}.

Denote by $S_{a,b}(t_1,t_2,\dots,t_b)$ the semihexagon of side-lengths $a,b,a+b,b$ (in counter-clock wise order, from the top side) with the dents at the positions $t_1,t_2,\dots,t_b$ on the base.
 By, say, \cite[equation (7.105)]{Stanley} and Cohn--Larsen--Propp's enumeration \cite[Proposition 2.1]{CLP}, we have\footnote{The notation $\textbf{1}^n$ in the argument of the Schur function  $s_{\lambda(\{t_1,\dots,t_b\})}(\textbf{1}^{b})$ stands for $n$ arguments equal to 1.}
\begin{equation}\label{CLPeq1}
\M(S_{a,b}(t_1,t_2,\dots,t_b))=\prod_{1\leq i< j \leq b}\frac{t_j-t_i}{j-i}=\textbf{s}_{\lambda(\{t_1,\dots,t_b\})}(\textbf{1}^{b}),
\end{equation}
 where  $\lambda(\{t_1,\dots,t_b\})$ is the the partition $(t_b-b+1,\dotsc,t_2-1,t_1)$. On the other hand, it is not hard to see that we also have
\begin{align}\label{CLPeq2}
\M\left(H_{x,y}(U;D)\right)&=\sum_{|S|=y}\M(S_{x+n-u,u+y}(U\cup S))\M(S_{x+n-d,d+y}(D\cup S))\notag\\
&=\sum_{|S|=y}s_{\lambda(U\cup S)}(\textbf{1}^{u+y})s_{\lambda(D \cup S)}(\textbf{1}^{d+y}),
\end{align}
where the sum runs over all $y$-subsets $S$ of  the complement of $U\cup D$.

 Indeed,  each tiling of $H_{x,y}(U;D)$ contains exactly $y$ vertical lozenges along the axis $l$. Moreover, these vertical lozenges must be at the positions in the complement of $U\cup D$.  Grouping those tilings which correspond to the same set of vertical lozenges, one could write the tiling number $\M(H_{x,y}(U;D))$ as the sum of tiling numbers:
\[\M(H_{x,y}(U;D))=\sum_{|S|=y}\M(H^{S}_{x,y}(U;D)),\]
where $H^{S}_{x,y}(U;D)$ denotes the region obtained from $H_{x,y}(U;D)$ by removing $y$ vertical lozenges at the positions in $S$. It is easy to see that each tiling of $H^{S}_{x,y}(U;D)$ can be separated into two tilings of two semi-hexagons by the axis $l$. It means that 
\[\M(H^{S}_{x,y}(U;D))=\M(S_{x+n-u,u+y}(U\cup S))\M(S_{x+n-d,d+y}(D\cup S)).\]
This implies the first identity in (\ref{CLPeq2}); the second identity follows from (\ref{CLPeq1}).

Applying (\ref{CLPeq2}), one can write our Shuffling Theorem in terms of Schur functions as
\begin{equation}\label{schureq}
\frac{  \sum_{|S|=y}\textbf{s}_{\lambda(U\cup S)}(\textbf{1}^{u+y})\textbf{s}_{\lambda(D \cup S)}(\textbf{1}^{d+y})}{  \sum_{|S|=y}\textbf{s}_{\lambda(U'\cup S)}(\textbf{1}^{u+y})\textbf{s}_{\lambda(D' \cup S)}(\textbf{1}^{d+y})}= \frac{\textbf{s}_{\lambda(U)}(\textbf{1}^{u})\textbf{s}_{\lambda(D)}(\textbf{1}^{d})}{\textbf{s}_{\lambda(U')}(\textbf{1}^{u})\textbf{s}_{\lambda(D')}(\textbf{1}^{d})},
  \end{equation}
where the sums are taken over all $y$-subsets $S$ of  the complement of $U\cup D$.

The $q$-analog of the Shuffling Theorem in \cite{shuffling} implies that (\ref{schureq}) still holds (up to a $q$-power) when $\textbf{1}^{n}$ is replaced by the sequence $(q,q^2,q^3,\dots,q^{n})$. It would be interesting to know if there is a more general Schur function identity behind (\ref{schureq}).

\begin{prob}\label{problem10}
Find a general Schur function identity behind the shuffling phenomenon. More particularly, find a `\emph{simple condition}' under which we have the following identity
\begin{equation}\label{schureq2}
\frac{  \sum_{|S|=y}\textbf{s}_{\lambda(U\cup S)}(\textbf{X}^{m+y})\textbf{s}_{\lambda(D \cup S)}(\textbf{X}^{m+(d-u)+y})}{  \sum_{|S|=y}\textbf{s}_{\lambda(U'\cup S)}(\textbf{X}^{m+y})\textbf{s}_{\lambda(D' \cup S)}(\textbf{X}^{m+(d-u)+y})}=C\cdot \frac{\textbf{s}_{\lambda(U)}(\textbf{X}^{m})\textbf{s}_{\lambda(D)}(\textbf{X}^{m+(d-u)})}{\textbf{s}_{\lambda(U')}(\textbf{X}^{m})\textbf{s}_{\lambda(D')}(\textbf{X}^{m+(d-u)})},
   \end{equation}
where $\textbf{X}^{m}$ denotes the sequence of variables $x_1,x_2,\dotsc,x_m$, and where $C$ is a monomial in $x_i$'s.
\end{prob}

It is worth noticing that Seok Hyun Byun \cite{Byun} and Markus Fulmek \cite{Ful} independently found a simple alternative proof for the Shuffling Theorem by using the connections between lozenge tilings and Schur function as in (\ref{CLPeq1}). While their proofs explain the identity (\ref{schureq}) and its $q$-analog,  they did not solve the problem above.

 By Cohn--Larson--Propp's enumeration again, our Shuffling Theorem  can be written purely in terms of tiling numbers as
\begin{align}\label{geointereq}
  \frac{\M(H_{x,y}(U; D))}{\M(H_{x,y}(U'; D'))}=\frac{\M(S_{x+y+n-u,u}(U))\M(S_{x+y+n-d,d}(D))}{\M(S_{x+y+n-u,u}(U'))\M(S_{x+y+n-d,d}(D'))}.
\end{align}
One readily sees that the two semi-hexagons in the numerator of the right-hand side of (\ref{geointereq}) are obtained by dividing the doubly-dented hexagon $H_{x+y,0}(U; D)$ along the horizontal axis $l$. More precisely, say by the Region-Splitting Lemma \cite{Tri18,Tri19}, we have
 \[\M(S_{x+y+n-u,u}(U))\M(S_{x+y+n-d,d}(D))=\M(H_{x+y,0}(U; D))\]
and
\[\M(S_{x+y+n-u,u}(U'))\M(S_{x+y+n-d,d}(D'))=\M(H_{x+y,0}(U'; D')).\]
 This means that identity  (\ref{geointereq}) could be rewritten as
\begin{equation}\label{geointereq2}
\M(H_{x,y}(U; D))\M(H_{x+y,0}(U'; D'))=\M(H_{x,y}(U'; D'))\M(H_{x+y,0}(U; D)).
\end{equation}
Both sides of (\ref{geointereq2}) count pairs of tilings of doubly-dented hexagons. There should be a bijection between these sets of pairs of tilings.

\begin{prob}\label{problem11}
Find a bijection that proves identity (\ref{geointereq2}).
\end{prob}

\begin{figure}\centering
\includegraphics[width=15cm]{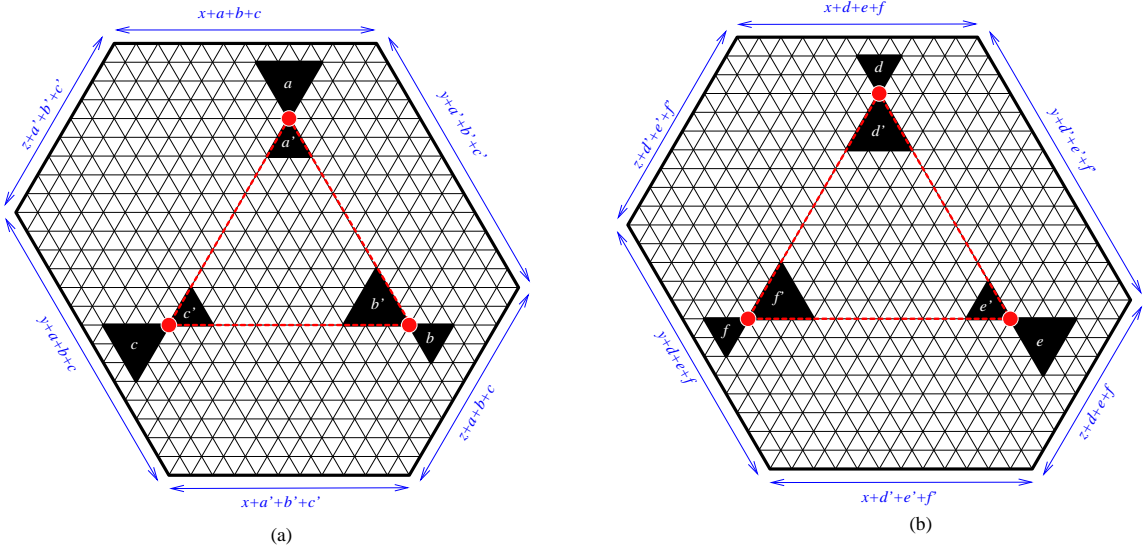}
\caption{Two sibling hexagons with a triad of bowties removed.}\label{Fig:triad}
\end{figure}

Several other examples of the shuffling phenomenon have been found. One of them is  the  shuffling phenomenon for  a hexagon with a removed ``\emph{triad of bowties}" \cite{CLR}. We consider the hexagon $H$ of side-lengths $x + a + b + c, y + a' + b' + c',z + a + b + c, x + a'+ b' + c', y + a + b + c, z + a' + b' + c'$. We remove from the interior of  $H$ three bowties shapes located at the vertices of an equilateral triangle $\Delta$ of side-length $t+a'+b'+c'$. Denote by $R=R^{\Delta}_{x,y,z}(a,a',b,b',c,c')$ the resulting region (see Figure \ref{Fig:triad}(a) for an example; the triangle with red dotted sides indicates the triangle $\Delta$).

 We now adjust  the side-lengths of the lobes in each bowtie so that the sum of the side-lengths is unchanged. Assume that the $(a,a')$-,  $(b,b')$-, and $(c,c')$-bowties are deformed into the $(d,d')$-, $(e,e')$-, and $(f,f')$-bowties, respectively, where $a+a'=d+d',$ $b+b'=e+e',$ and $c+c'=f+f'$. We now create a `sibling' region $R'$ of $R$ as follows. We start with the  hexagon $H'$ of side-lengths $x + d +e + f, y + d' + e' + f',z + d+e+f, x + d'+ e' + f', y + d + e + f, z + d' + e' + f'.$ We now remove the above three new bowties from $H'$  so that they are located at the vertices of a triangle $\Delta'$ of side-length $t+d'+e'+f'$ and that the distances from the  top, the left, and the right bowties to the north, the southeast, and the southwest sides of $H'$ are equal to the corresponding distances in the original region.  Denote by $R'=R^{\Delta'}_{x,y,z}(d,d',e,e',f,f')$ the new region.
 
It has been shown that the ratio of tiling  numbers of the above two regions is given by a simple product formula \cite{CLR}. We would also like to emphasize that, in general, the tiling numbers of the regions $R^{\Delta}_{x,y,z}(a,a',b,b',c,c')$ and $R^{\Delta'}_{x,y,z}(d,d',e,e',f,f')$ are \emph{not} given by simple product formulas (see Theorem 1 in \cite{CLR}).  See Figure \ref{Fig:triad} for an example of $x=4,y=6,z=2, t=4,a=3,a'=2,b=2,b'=3,c=3,c'=2,d=2,d'=3,e=3,e'=2,f=2,f'=3$; the triangles $\Delta$ and $\Delta'$ are the ones with red vertices in pictures (a) and (b), respectively. We note that it actually needs three more parameters to define the triangle $\Delta$ and that the triangle $\Delta'$ is uniquely determined by $\Delta$. It means that each of the regions $R$ and $R'$ depends on $12$ parameters.

We have realized that the above instance of the  shuffling phenomenon could be extended to the weighted case. In particular, we consider the generating functions with the elliptic weight $\wt_1=\frac{Xq^{j}+Yq^{-j}}{2}$ of $R^{\Delta}_{x,y,z}(a,a',b,b',c,c')$ and $R^{\Delta'}_{x,y,z}(d,d',e,e',f,f')$. Like the unweighted case, these tilling generating functions are not given by any simple product formulas; however, their ratio seems to be a simple product. 

\begin{prob}\label{problem12}
Prove that the ratio of tiling generating functions \[\frac{\M_1(R^{\Delta}_{x,y,z}(a,a',b,b',c,c'))}{\M_1(R^{\Delta'}_{x,y,z}(d,d',e,e',f,f'))}\] is always given by a simple product formula. Here, we use the notations $\M_i(R)$ for the tiling generating function of $R$ using the weight $\wt_i$, $i=1,2,3$, as described in the previous section.
\end{prob}

We want to point out that each tiling generating function in the above problem depends on \emph{fifteen} parameters. Proving this theorem would be technically challenging. We note that the conjectural formula for the above ratio has been found in the particular case when $X=2$ and $Y=0$ (i.e., the case of the natural weight $\wt_2$). However, a formula in the general case is still unknown.

Like its unweighted version, the conjecture would imply a number of results in the weighted enumeration of lozenge tilings. Let us point out one of them.  Rosengren in \cite{Rosen} found an intriguing consequence of his weighted enumeration of the cored hexagons (see Corollary 2.2 therein). Intuitively, if we reflect the  removed triangle through the center of the hexagon, then the tiling generating function of the core hexagon changes  by only a simple multiplicative factor. The conjecture (if proved) would give a conceptual explanation for his observation. Indeed, we consider the special case when $X=Y=1$, $a=d'$,  $a'=b=b'=c=c'=d=e=e'=f=f'=0$, and $q$ is replaced by $\sqrt{q}$. Then the ratio between the tiling generating functions of $R_{x,y,z}(a,0,0,0,0,0)$ and (the horizontal refection of) $R_{x,y,z}(0,a,0,0,0,0)$ is exactly the ratio in Rosengren's Corollary 2.2.

It has been shown that the original shuffling phenomenon also holds for the reflectively symmetric tilings and centrally symmetric tilings  of the doubly--dented hexagons \cite{shuffling2, shuffling3}. It suggests we do the same for hexagons with a triad of bowtie holes. In particular, we first focus on the \emph{cyclically symmetric tilings}, i.e., the tilings invariant under $120^{\circ}$ rotations. We note that the region $R^{\Delta}_{x,y,z}(a,a',b,b',c,c')$ must be cyclically symmetric itself in order to have a cyclically symmetric tiling. More precisely, we must have $x=y=z$, $a=b=c$, $a'=b'=c'$, and the triangle $\Delta$ must be at the center of the region. To emphasize the symmetric location of the triangle, we use the notation $\Delta_0$ instead of $\Delta$.

\begin{prob}\label{problem13}
Let $x,a,a',d,d'$ be non-negative integers so that $a+a'=d+d'$. Find a formula for the ratio of tiling generating functions \[\frac{\M_c(R^{\Delta_0}_{x,x,x}(a,a',a,a',a,a'))}{\M_c(R^{\Delta'_0}_{x,x,x}(d,d',d,d',d,d'))}.\] Here we use the notation $\M_c(R)$ for the weighted sum of cyclically symmetric  tilings  of $R$ using the weight $\wt_1=\frac{Xq^{j}+Yq^{-j}}{2}$ (or some specialization).
\end{prob}

We are also interested in the shuffling phenomenon for the reflectively symmetric tilings of hexagons with a triad of bowties.  The region $R^{\Delta}_{x,y,z}(a,a',b,b',c,c')$ must be reflectively symmetric itself in order to have a reflectively symmetric tiling. More precisely, we must have $y=z$, $b=c$, $b'=c'$, and the triangle $\Delta$ must be on the symmetry axis. To emphasize the symmetric location of the triangle, we use the notation $\Delta_1$ instead of $\Delta$. 

\begin{prob}\label{problem14}
Let $x,y,a,a',b,b',d,d',e,e'$ be non-negative integers so that $a+a'=d+d'$ and $b+b'=e+e'$.  Find a formula for the ratio of tiling generating functions \[\frac{\M_r(R^{\Delta_1}_{x,y,y}(a,a',b,b',b,b'))}{\M_r(R^{\Delta'_1}_{x,y,y}(d,d',e,e',e,e'))}.\] Here we use the notation $\M_r(R)$ for the weighted sum of reflectively symmetric  tilings  of $R$ using the weight $\wt_3=\frac{q^{j}+q^{-j}}{2}$.
\end{prob}

We note that in Problem \ref{problem14}, the tilings is weighted by $\wt_3=\frac{q^{j}+q^{-j}}{2}$, a special case of the weight $\wt_1=\frac{Xq^{j}+Yq^{-j}}{2}$ in Problem \ref{problem13}.  Our data suggests that  the general weight $\wt_1$ does not yield a nice product formula for the tiling ratio in this case.

\begin{figure}\centering
\includegraphics[width=12cm]{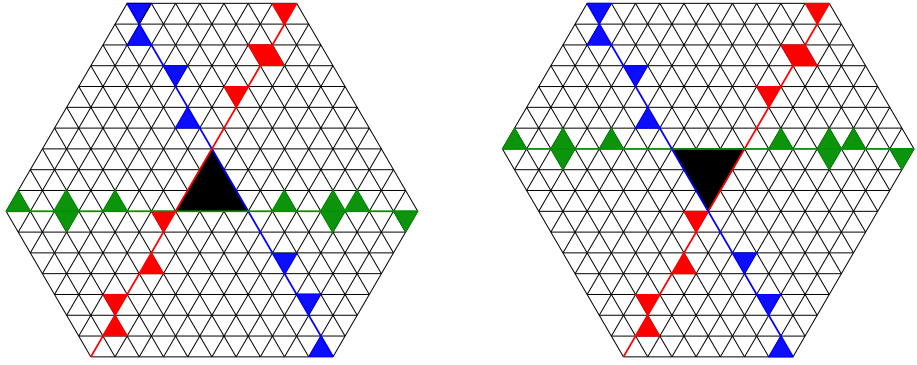}
\caption{The hexagons with holes on three crossing lines.}\label{Fig.threeline}
\end{figure}

Byun generalized the shuffling phenomenon for a triad of bowties in \cite{CLR}. He showed that the tiling number of a hexagon with holes on three crossing lines only changes by a simple multiplicative factor if we flip the central triangular hole and then translate these lines of holes (see Theorem 2.1 in \cite{Byun2}). See Figure \ref{Fig.threeline} for an illustration. More precisely, Byun's region is obtained from a cored-hexagon of side-lengths $n,n+x,n,n+x,n,n+x$ (in clockwise order, starting from the north side) with the (up-pointing) triangular hole of side-length $x$ in the center. We next extend the sides of the central hole to three equal segments distinguished by three different colors in the left region. On these segments, we put some additional holes. Byun uses \emph{twelve} sequences $A_1,A_2,\dots,A_6,B_1,B_2,\dots,B_6$ to record the position of the additional holes, for brevity, we denote $\textbf{A}:=(A_1,\dots,A_6)$ and $\textbf{B}:=(B_1,\dots,B_6)$. The resulting region is denoted by $H_{n,x}(\textbf{A},\textbf{B})$. We now consider the hexagon of side-lengths $n+x,n,n+x,n,n+x,n$ (also in clockwise order, starting from the north side) with a \emph{down-pointing} central hole of size $x$. We now implant the three segments of holes in the previous region on the extended sides of the new central hole in the new core-hexagon. Denote the new region by $\overline{H}_{n,x}(\textbf{A},\textbf{B})$ (see the right region in the figure). 

One would like to generalize this exciting result to the tiling generating function with respect to the weight $\wt_1$. More precise, we assign weight to lozenges in the hexagon with holes on three crossing lines $H_{n,x}(\textbf{A},\textbf{B})$ (resp., $\overline{H}_{n,x}(\textbf{A},\textbf{B})$) using the weight $\wt_1$ (or some specialization of its), with the horizontal axis running along the base and the vertical axis passing the top (resp., bottom) vertex of the central triangular hole.  We conjecture that the ratio of the two resulting generating functions is given by a simple product formula.

\begin{prob}\label{problem15}
Generalizing Byun's Theorem 2.1(a) in \cite{Byun2} to tiling generating functions using the weight $\wt_1$ (or some specialization).
\end{prob}

Inspired by Problems \ref{problem13} and \ref{problem14} above, we are also interested in the symmetric version of Problem \ref{problem15} above.

\begin{prob}\label{problem16}
(a) Generalizing Byun's Theorem 2.1(b) in \cite{Byun2} to generating functions of cyclically symmetric tilings using the weight $\wt_1$ (or some specialization).

(b) Investigate the ratio of generating functions of reflectively symmetric tilings of hexagons with holes on three crossing lines.
\end{prob}

As shown in \cite{Tri1, Tri6}, in many cases, one can convert a lozenge tiling problem to a domino tiling problem, and vice versa. It suggests the existence of a shuffling theorem for domino tilings.

\begin{prob}\label{problem17}
Find an example of the shuffling phenomenon for domino tilings.
\end{prob}

\section{Connections to Electrical networks}
This section is devoted to connections between tilings and other mathematical areas. Among many interesting connections, we focus on the connection to the study of electrical networks. Again, we recommend the reader the excellent survey paper of James Propp \cite{Propp2} for various connections and applications of the enumeration of tilings.

The study of \emph{electrical networks} comes from classical physics with the work of Ohm and Kirchhoff more than 100 years ago.
A \emph{circular planar electrical network} (or simply \emph{electrical network}) is a graph $G=(E,V)$ embedded in a disk with a set of distinguished vertices $N\subseteq V$ on the circle, called \emph{nodes}, and a \emph{conductance} function $wt: E\rightarrow \mathbb{R}^+$. The electrical networks were first studied systematically by Colin de Verdi\`{e}re \cite{Colin1} and Curtis, Ingerman, Moores, and Morrow \cite{Curtis1,Curtis2}.  We refer the reader to, e.g., for the recent development of the topic \cite{Alman,HK,KW, Lam2,LP1,Yi14}.

\begin{figure}\centering
\includegraphics[width=10cm]{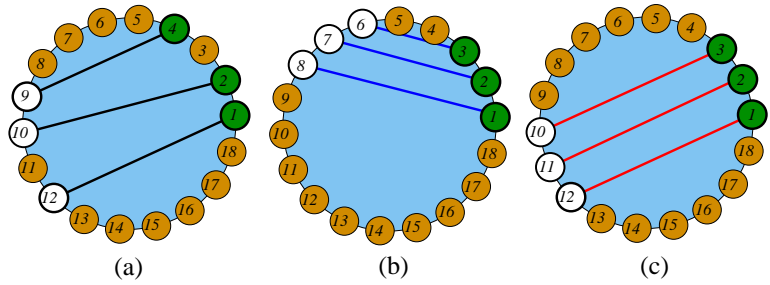}
\caption{(a) A non-contiguous minor.  (b) A contiguous minor that is not a central minor. (c) A central minor.}
\label{Fig:minor}
\end{figure}

We arrange the indices $1,2,\dotsc, n$ of a $n\times n$ matrix $M$ in counter-clockwise order around the circle.  Let $A=\{a_1,a_2,\dotsc,a_k\}$ and $B=\{b_1,b_2,\dotsc,b_k\}$ be two sets of indices so that $a_1,$ $a_2,$ $\dotsc,$ $a_k,$ $b_k,$ $b_{k-1},$ $\dotsc,$ $b_1$ appear in counter-clockwise order around the circle. The \emph{circular minor} $\det M_{A}^{B}$ is defined to be the minor of $M$ obtained from the rows $a_1,a_2,\dotsc,a_k$ and the columns $b_k,b_{k-1},\dots,b_1$. When $A$ and $B$ are non-interlaced around the circle, we can represent the circle minor $\det M_{A}^{B}$ by a disk diagram with $k$ chords connecting notes $a_i$ to $b_i$. See Figure \ref{Fig:minor} for examples.

\begin{figure}\centering
\includegraphics[width=13cm]{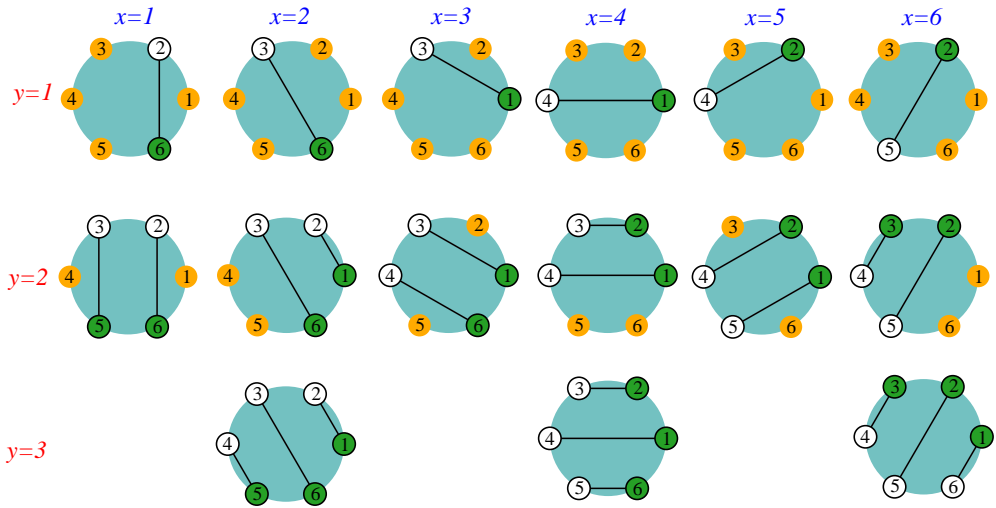}
\caption{All 15 small central minors when $n=6$.}\label{Fig:smallcentral}
\end{figure}

A \emph{contiguous minor} of  a matrix $M$ is a circular minor whose row indices and whose column indices are contiguous on the circle. A \emph{(small) central minor} is a non-interlaced contiguous minor whose row indices and column indices are opposite (or almost opposite depending on the parity of $n$) around the circle. There are $\binom{n}{2}$ central minors, whether $n$ is even or odd. See Figures \ref{Fig:minor} and \ref{Fig:smallcentral} for examples of these special types of minors. 

Associated with an electrical network  is a \emph{response matrix} that measures the response of the network to potential applied at the nodes.
It has been shown that a matrix $M$ is the response matrix of an electrical network if and only if it is symmetric with row and column sums equal to zero and each circular minor $\det M_{A}^{B}$ is non-negative (see Theorem 4 in \cite{Curtis1}). 
An electrical network is called \emph{well-connected} if, for any two non-interlaced sets of nodes $A$ and $B$, there are  $k$ pairwise vertex-disjoint paths in $G$ connecting the nodes in $A$ to the nodes in $B$, where $|A|=|B|=k$.

\begin{figure}\centering
\includegraphics[width=6cm]{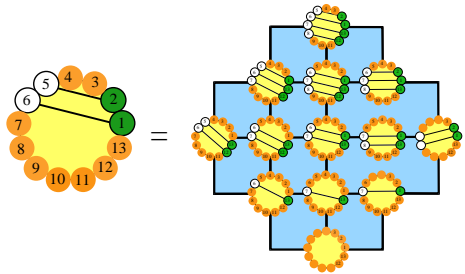}
\caption{Encoding a contiguous minor as a weighted Aztec diamond.}
\label{Fig:KWfig}
\end{figure}

 R. Kenyon and D. Wilson \cite{KW} generalize the work of  Colin de Verdi\`{e}re \cite{Colin1}  by showing how to test the well-connectivity of an electrical network  by checking the positivity of the $\binom{n}{2}$ central minors of the response matrix.
The test is based on their interesting finding of the connection between linear algebra, electrical networks, and tilings.
\begin{thm}[Kenyon--Wilson's Theorem \cite{KW}]
Any contiguous minor can be written as a Laurent polynomial in central minors. Moreover, this Laurent polynomial is the generating function of tilings of a (weighted) `\emph{Aztec diamond}.' 
\end{thm}
 See Figure \ref{Fig:KWfig} for an illustration of the theorem.

The \emph{Aztec diamond} of order $h$ with the center located at the lattice point $(x_0,y_0)$ is the region consisting of all unit squares inside the contour $|x-x_0|+|y-y_0|\leq h+1$. Elkies, Kuperberg, Larsen, and Propp \cite{Elkies1, Elkies2} have shown that the number of domino tilings of the Aztec diamond of order $h$ is exactly $2^{h(h+1)/2}$. This work has inspired a large body of work in the enumeration of tilings.

\begin{figure}\centering
\includegraphics[width=13cm]{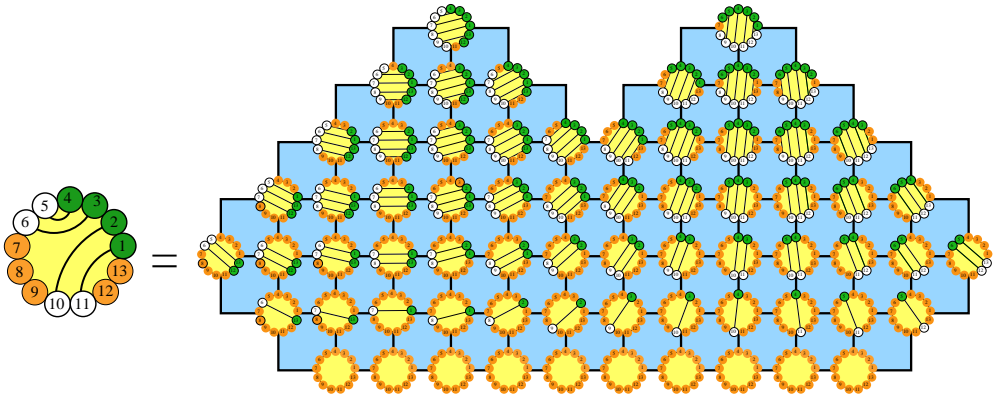}
\caption{Illustration of the correspondence between a semicontiguous minor (left-hand side) and the domino tilings of a region weighted by central minors (right-hand side). The lattice points on the right-hand side are replaced by the corresponding central minors. This picture first appeared in \cite{Tri20}.}
\label{CMgraphic1}
\end{figure}

We are also interested in a larger family of minors, called \emph{semicontiguous minors}. A \emph{semicontiguous minor} is a circular minor $\det M_{A}^{B}$, where \emph{at least one} of $A$ and $B$ is contiguous. Kenyon and Wilson conjectured that any semicontiguous minor could also be written as the tiling generating function of some region on the square lattice. This conjecture was recently proved in \cite{Tri20} by building a special family of regions on the square lattice whose tiling generating functions are given by the semicontiguous minors. In particular, our region is obtained from an \emph{Aztec rectangle} (a natural generalization of the Aztec diamond) or a special union of two Aztec diamonds by trimming the base along a zigzag path determined by the non-contiguous index set. It would be interesting to know if any general circular minor can be encoded as domino tilings of  a region.

\begin{prob}\label{problem18} 
Can any circular minor be written as the tiling generating function of a region on a square lattice? If the answer is ``NO," characterize all such circular minors.
\end{prob}

\section{Other problems}

We consider the number of cyclically symmetric tilings of two families of hexagons with four triangles removed as follows. 

\begin{figure}\centering
\includegraphics[width=13cm]{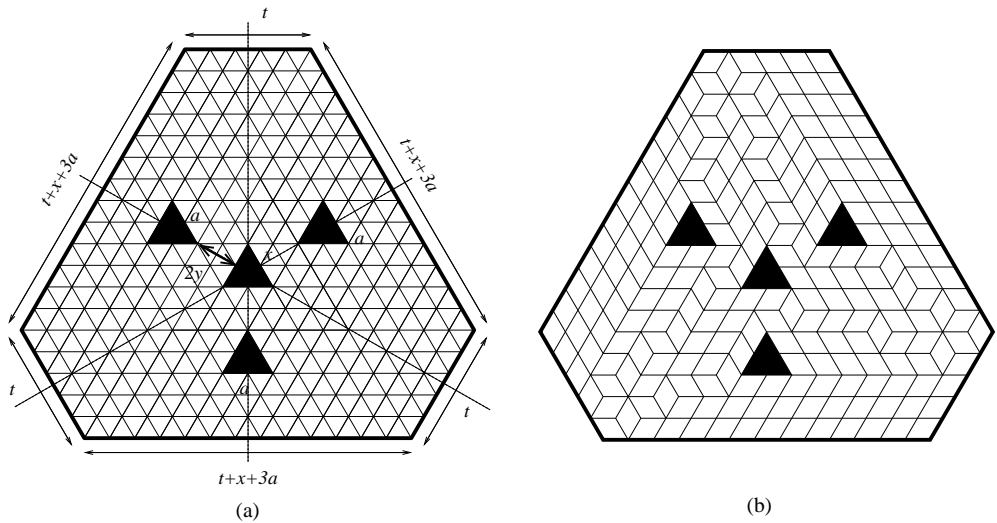}
\caption{(a) The hexagon with four holes $\mathcal{H}_{5,1}(2,2)$. (b) A cyclically symmetric tiling of $\mathcal{H}_{5,1}(2,2)$. This picture first appeared in \cite{LR2}.}\label{4holes}
\end{figure}

Let $x,y,t,a$ be non-negative integers. Our first family consists of the hexagons with side-lengths $t+x+3a,\ t,\ t+x+3a,\ t,\ t+x+3a,\ t$,  in which an up-pointing triangle of side-length $x$ has been removed from the center. In addition, three up-pointing triangles of side-length $a$ have been removed in a symmetric way along the intervals connecting the center to the midpoints of the southern, northeastern, and northwestern sides of the hexagon. We assume besides that the distance from the central hole to each of the three satellite holes is $2y$. Denote by $\mathcal{H}_{t,y}(a,x)$ the resulting region (see Figure \ref{4holes} for an example; the black triangles represent the triangles that have been removed).

\begin{figure}\centering
\includegraphics[width=13cm]{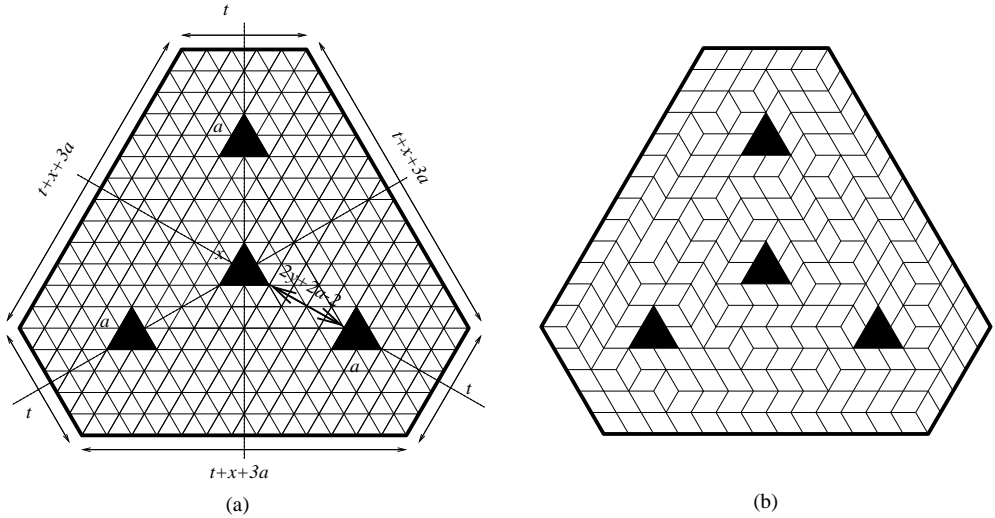}
\caption{(a) The hexagon with four holes $\overline{\mathcal{H}}_{5,1}(2,2)$. (b) A cyclically symmetric tiling of $\overline{\mathcal{H}}_{5,1}(2,2)$. This picture first appeared in \cite{LR2}.}\label{4holes2} 
\end{figure}

The second family also consists of hexagons with four similar triangular holes; however, the $a$-triangles now lie on the other side of the center, as shown in Figure \ref{4holes2}. The resulting region is denoted by $\overline{\mathcal{H}}_{t,y}(a,x)$. Two simple product formulas for the numbers of  cyclically symmetric tilings of $\mathcal{H}_{t,y}(a,x)$ and $\overline{\mathcal{H}}_{t,y}(a,x)$ were provided in   \cite{LR2} .

Recently, M. Ciucu and I. Fischer \cite{CF19} conjecture a striking connection between the whole number of tilings and the number of cyclically symmetric tilings of $\mathcal{H}_{t,y}(a,x)$.

\begin{prob}[Conjecture 1 in \cite{CF19}] \label{problem19}
Show that
\begin{align}\label{CFeq}
\frac{\M(\mathcal{H}_{t,y}(a,x))}{\M_r(\mathcal{H}_{t,y}(a,x))^3}=\frac{\M(\mathcal{H}_{t,0}(a,x))}{\M_r(\mathcal{H}_{t,0}(a,x))^3}\left[\prod_{i=1}^{y}\frac{(x+6i-4)(x+3a+6i-2)}{(x+6i-2)(x+3a+6i-4)}\right]^2,
\end{align}
where $\M_r(R)$ denotes the number of cyclically symmetric tilings of $R$.
\end{prob}
We note that when $y=0$, all four holes in $\mathcal{H}_{t,0}(a,x)$ are glued together, and  the region has the same tiling number as a cored hexagon. Recall that the number of tilings of a cored hexagon is given by  a simple product formula in \cite{CEKZ}, so the expression on the right-hand side of (\ref{CFeq}) can be expressed as a simple product (see the discussion before Conjecture 2 in \cite{CF19}). As the number of cyclically symmetric tilings of $\mathcal{H}_{t,y}(a,x)$ has been found, one only needs to find the tiling number of $\mathcal{H}_{t,y}(a,x)$ to prove the conjecture. However,  as discussed in \cite{CF19}, this task would not be easy.

The similarity of $\mathcal{H}_{t,y}(a,x)$ and $\overline{\mathcal{H}}_{t,y}(a,x)$ suggests the existence of a nice formula for the ratio $\frac{\M(\overline{\mathcal{H}}_{t,y}(a,x))}{\M_c(\overline{\mathcal{H}}_{t,y}(a,x))^3}$. However, such a formula was not provided in \cite{CF19}. It would be interesting to obtain that formula.

\begin{prob}\label{problem20}
Find a formula for the ratio of tiling numbers
\begin{align}\label{CFeq2}
\frac{\M(\overline{\mathcal{H}}_{t,y}(a,x))}{\M_c(\overline{\mathcal{H}}_{t,y}(a,x))^3}.
\end{align}
\end{prob}

It is worth noticing that when $a=x=0$, i.e., when all four triangular holes are all vanished, the equation (\ref{CFeq}) becomes
\[\frac{\M(Hex(n,n,n))}{\M_c(Hex(n,n,n))^3}=\left[\frac{\left(\frac{1}{3}\right)_n}{\left(\frac{2}{3}\right)_n}\right]^2,\]
 which follows from the well-known enumerations of symmetry classes of plane partitions (see, i.e., \cite{Stanley2}). Here we use the Pochhammer symbol $(x)_n:=x(x+1)(x+2)\cdots(x+n-1)$.
In particular, the tiling number in the numerator is precisely the number of boxed plane partitions;  the tiling number in the denominator is the number of cyclically symmetric plane partitions. Macdonald conjectured the following weighted enumeration of cyclically symmetric plane partitions \cite[Ex. 18, p. 85]{Macdonald}:
\begin{equation}\label{Macdonaleq}
\sum_{\pi}q^{|\pi|}=\prod_{i=1}^{n}\frac{1-q^{3i-1}}{1-q^{3i-2}}\prod_{1\leq i<j\leq n}\frac{1-q^{3(2i+j-1)}}{1-q^{3(2i+j-2)}}\prod_{1\leq i<j,k\leq n}\frac{1-q^{3(i+j+k-1)}}{1-q^{3(i+j+k-2)}},
\end{equation}
where the sum is over all cyclically symmetric plane partitions $\pi$ that are contained in a $(n\times n\times n)$-box. The unweighted version of the conjecture was proved by Andrews \cite{Andrews}; the full conjecture was proved by Mills, Robbins, and Rumsey \cite{MRR}.

\begin{figure}\centering
\includegraphics[width=7cm]{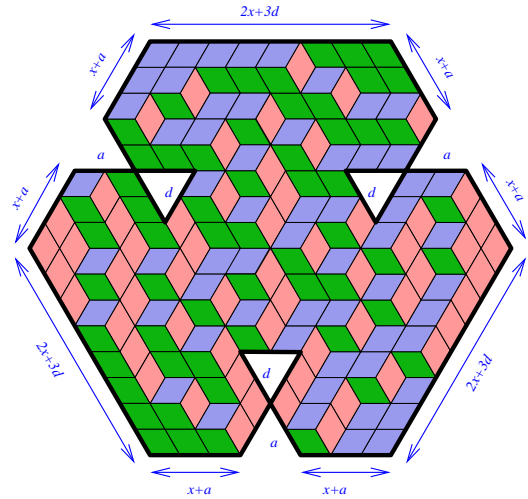}
\caption{A cyclically symmetric tiling of the region with three bowties removed.}\label{Fig:symmetric3bowtie}
\end{figure}

At this point, one would think about similar weighted versions of the tiling numbers in (\ref{CFeq}). Unfortunately, there have been no such nice $q$-enumerations in general. However, in the case when $x = 0$ (i.e., the central hole vanishes) and the three satellite holes are attached to the hexagon's boundary, we have a nice $q$-enumeration, as claimed in Theorem \ref{Bowtiethm}.  We want to obtain a formula for the weighted sum of  the cyclically symmetric tilings of the region $F\begin{pmatrix} x&x&x\\a&a&a\\d&d&d\end{pmatrix}$ (an $F$-type region in Theorem  \ref{Bowtiethm} must have $a=b=c$, $d=e=f$, $x=y=z$ in order to have  a cyclically symmetric tiling). See Figure \ref{Fig:symmetric3bowtie} for an example of such tiling. Equivalently, we are interested in finding  the following volume generating function:

\begin{prob}\label{problem21}
Find the volume generating function of the cyclically symmetric stacks fitting in the compound room $\mathcal{B}\begin{pmatrix} x&x&x\\a&a&a\\d&d&d\end{pmatrix}$.
\end{prob}

\begin{figure}\centering
\includegraphics[width=15cm]{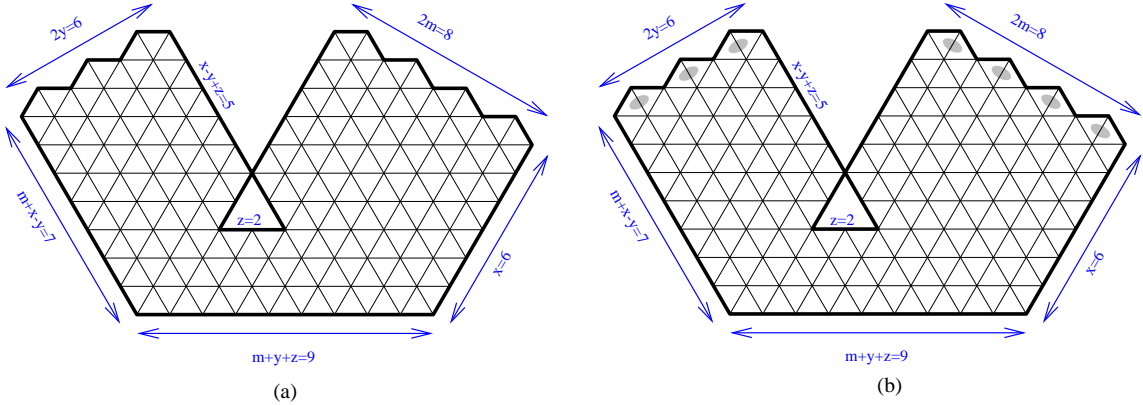}
\caption{Two regions enumerated by Ciucu and Fischer: (a) $D_{6,3,2,4}$ and (b) $D'_{6,3,2,4}$. The lozenges with shaded cores are weighted by $1/2$.}\label{Fig:CFcutoff}
\end{figure}

\medskip

Ciucu and Fischer \cite[Theorems 2.3 and 2.4]{CF15} proved and generalized two conjectures of Ciucu and Krattenthaler \cite[Conjectures A.1 and A.2]{CK02} by enumerating the regions $D_{x,y,z,m}$ and $D'_{x,y,z,m}$ in Figure \ref{Fig:CFcutoff}.  These regions look essentially the same, and the only difference is that the lozenges running along two zigzag cuts in $D'_{x,y,z,m}$ are weighted by $1/2$ (see the lozenges with shaded cores). One can view these regions as one-third of an $F$-type region in the previous problem. From the tiling formulas, we could realize that the weighted enumeration of tilings of $D'_{x,y,z,m}$ is obtained from the tiling formula of $D_{x,y,z,m}$ by replacing $x$ by $x-1/2$. This property reminds us of the ``\emph{combinatorial reciprocity phenomenon}": even though the region $D_{x,y,z,m}$ is \emph{not} defined when $x$ is a half-integer, its tiling formula gives the number of combinatorial objects of a different sort (here are the tilings of $D'_{x,y,z,m}$) when evaluated at half-integers. We refer the reader to, e.g., \cite{Beck,StanleyRecip,ProppRecip} for more discussions about the phenomenon. 

\begin{prob}\label{problem22}
Explain combinatorially the above  $\left(\frac{1}{2}\right)$-phenomenon for the tiling numbers of the region $D_{x,y,z,m}$ and $D'_{x,y,z,m}$.
\end{prob}

\begin{figure}\centering
\setlength{\unitlength}{3947sp}%
\begingroup\makeatletter\ifx\SetFigFont\undefined%
\gdef\SetFigFont#1#2#3#4#5{%
  \reset@font\fontsize{#1}{#2pt}%
  \fontfamily{#3}\fontseries{#4}\fontshape{#5}%
  \selectfont}%
\fi\endgroup%
\resizebox{!}{12cm}{
\begin{picture}(0,0)%
\includegraphics{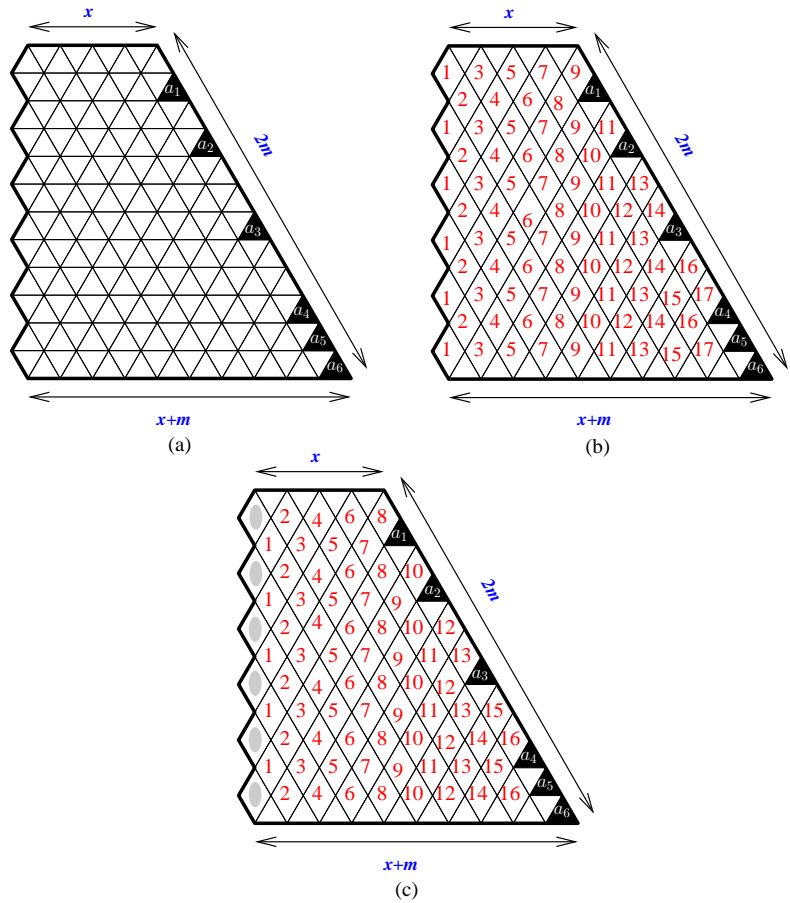}%
\end{picture}%
%
%

\begin{picture}(9893,11549)(775,-10942)
\put(7261,-9106){\makebox(0,0)[lb]{\smash{{\SetFigFont{14}{16.8}{\rmdefault}{\bfdefault}{\updefault}{\color[rgb]{1,1,1}$a_4$}%
}}}}
\put(3174,-1314){\makebox(0,0)[lb]{\smash{{\SetFigFont{14}{16.8}{\rmdefault}{\bfdefault}{\updefault}{\color[rgb]{1,1,1}$a_2$}%
}}}}
\put(3774,-2379){\makebox(0,0)[lb]{\smash{{\SetFigFont{14}{16.8}{\rmdefault}{\bfdefault}{\updefault}{\color[rgb]{1,1,1}$a_3$}%
}}}}
\put(4389,-3429){\makebox(0,0)[lb]{\smash{{\SetFigFont{14}{16.8}{\rmdefault}{\bfdefault}{\updefault}{\color[rgb]{1,1,1}$a_4$}%
}}}}
\put(4599,-3781){\makebox(0,0)[lb]{\smash{{\SetFigFont{14}{16.8}{\rmdefault}{\bfdefault}{\updefault}{\color[rgb]{1,1,1}$a_5$}%
}}}}
\put(4809,-4126){\makebox(0,0)[lb]{\smash{{\SetFigFont{14}{16.8}{\rmdefault}{\bfdefault}{\updefault}{\color[rgb]{1,1,1}$a_6$}%
}}}}
\put(7471,-9458){\makebox(0,0)[lb]{\smash{{\SetFigFont{14}{16.8}{\rmdefault}{\bfdefault}{\updefault}{\color[rgb]{1,1,1}$a_5$}%
}}}}
\put(7681,-9803){\makebox(0,0)[lb]{\smash{{\SetFigFont{14}{16.8}{\rmdefault}{\bfdefault}{\updefault}{\color[rgb]{1,1,1}$a_6$}%
}}}}
\put(5626,-6278){\makebox(0,0)[lb]{\smash{{\SetFigFont{14}{16.8}{\rmdefault}{\bfdefault}{\updefault}{\color[rgb]{1,1,1}$a_1$}%
}}}}
\put(6046,-6991){\makebox(0,0)[lb]{\smash{{\SetFigFont{14}{16.8}{\rmdefault}{\bfdefault}{\updefault}{\color[rgb]{1,1,1}$a_2$}%
}}}}
\put(8080,-609){\makebox(0,0)[lb]{\smash{{\SetFigFont{14}{16.8}{\rmdefault}{\bfdefault}{\updefault}{\color[rgb]{1,1,1}$a_1$}%
}}}}
\put(8500,-1322){\makebox(0,0)[lb]{\smash{{\SetFigFont{14}{16.8}{\rmdefault}{\bfdefault}{\updefault}{\color[rgb]{1,1,1}$a_2$}%
}}}}
\put(9100,-2387){\makebox(0,0)[lb]{\smash{{\SetFigFont{14}{16.8}{\rmdefault}{\bfdefault}{\updefault}{\color[rgb]{1,1,1}$a_3$}%
}}}}
\put(9715,-3437){\makebox(0,0)[lb]{\smash{{\SetFigFont{14}{16.8}{\rmdefault}{\bfdefault}{\updefault}{\color[rgb]{1,1,1}$a_4$}%
}}}}
\put(9925,-3789){\makebox(0,0)[lb]{\smash{{\SetFigFont{14}{16.8}{\rmdefault}{\bfdefault}{\updefault}{\color[rgb]{1,1,1}$a_5$}%
}}}}
\put(10135,-4134){\makebox(0,0)[lb]{\smash{{\SetFigFont{14}{16.8}{\rmdefault}{\bfdefault}{\updefault}{\color[rgb]{1,1,1}$a_6$}%
}}}}
\put(6646,-8056){\makebox(0,0)[lb]{\smash{{\SetFigFont{14}{16.8}{\rmdefault}{\bfdefault}{\updefault}{\color[rgb]{1,1,1}$a_3$}%
}}}}
\put(2754,-601){\makebox(0,0)[lb]{\smash{{\SetFigFont{14}{16.8}{\rmdefault}{\bfdefault}{\updefault}{\color[rgb]{1,1,1}$a_1$}%
}}}}
\end{picture}}
\caption{(a) The quartered hexagon with dents on the right side. (b) How to assign weights to lozenges in $Q_{4}(2,4,7,10,11,12)$. (c) How to assign weights to lozenges in $Q'_{4}(2,4,5,10,11,12)$.}\label{Fig:Semitwodent2}
\end{figure}

A similar thing happens in \cite[Theorems 1.3 and 1.4]{Semitwodent}. In Theorem 1.3, we consider the ratio of tiling generating functions of pairs of  (dented) quartered hexagons  $Q_{x}(\textbf{a})$ and $Q_{y}(\textbf{a})$. They have the same dents on the right sides at the positions in $\textbf{a}=(a_i)_{i=1}^{m}$; the only difference is at their widths. We assign weights to lozenges of  $Q_{x}(\textbf{a})$ and $Q_{y}(\textbf{a})$ as in Figure \ref{Fig:Semitwodent2}(b): each vertical lozenges with label $x$ is weight by $(q^{x}+q^{-x})/2$, for $x\geq 1$ (lozenges of other orientations are all weighted by $1$). The pair of quartered hexagons  $Q'_{x}(\textbf{a})$ and $Q'_{y}(\textbf{a})$ in  Theorem 1.4 are the same as  that in Theorem 1.3; the only difference is the weights of lozenges. Figure \ref{Fig:Semitwodent2}(c) illustrates the weight assignment of lozenges in  $Q'_{x}(\textbf{a})$ and $Q'_{y}(\textbf{a})$ (the lozenges with shaded cores on the left side are weighted by $1/2$). 
We have
\begin{align}\label{Qeq1}
\frac{\M(Q_{x}((a_i)_{i=1}^{m}))}{\M(Q_{y}((a_i)_{i=1}^{m}))}&=q^{2(y-x)(\sum_{i=1}^{m}a_i- m^2)}\prod_{i=1}^{m}\frac{(q^{2(2y+a_i+1)};q^2)_{2i-a_i-1}}{(q^{2(2x+a_i+1)};q^2)_{2i-a_i-1}}
\end{align}
and
\begin{align}\label{Qeq2}
\frac{\M(Q'_{x}((a_i)_{i=1}^{m}))}{\M(Q'_{y}((a_i)_{i=1}^{m}))}&=q^{2(y-x)(\sum_{i=1}^{m}a_i- m^2)}\prod_{i=1}^{m}\frac{(q^{2(2y+a_i)};q^2)_{2i-a_i-1}}{(q^{2(2x+a_i)};q^2)_{2i-a_i-1}}.
\end{align}
It is easy to see that $\frac{\M(Q'_{x}((a_i)_{i=1}^{m}))}{\M(Q'_{y}((a_i)_{i=1}^{m}))}$ is obtained from $\frac{\M(Q_{x}((a_i)_{i=1}^{m}))}{\M(Q_{y}((a_i)_{i=1}^{m}))}$ by replacing $x$ by $x-1/2$ and $y$ by $y-1/2$. It would be interesting to find a direct explanation for this, i.e., we want an explanation without requiring any calculation of the ratios tiling generating functions 

\begin{prob}\label{problem23}
Explain the $\left(\frac{1}{2}\right)$-phenomenon for the ratios of tiling generating functions of the quartered hexagons $Q_{x}((a_i)_{i=1}^{m})$ and $Q'_{x}((a_i)_{i=1}^{m})$.
\end{prob}

\begin{figure}\centering
\includegraphics[width=6cm]{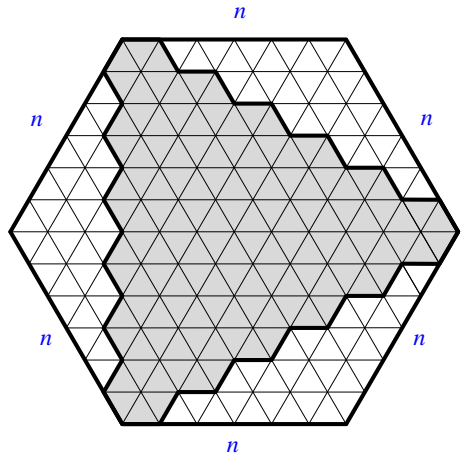}
\caption{The region $\mathcal{TT}_6$.}\label{Fig:Tn}
\end{figure}

In their excellent paper about the tilling enumeration of  hexagons with a maximal corner cut off \cite{CK02}, Ciucu and Krattenthaler found an unusual pattern for the tiling number of a triangular region denoted by $\mathcal{TT}_n$ (see the shaded region in Figure \ref{Fig:Tn} for $\mathcal{TT}_6$). The number of  tilings of $\mathcal{TT}_n$ factors as follows for $n\leq  7$:
\begin{align}
&\M(\mathcal{TT}_1)=2;\notag\\
&\M(\mathcal{TT}_2)=3^2;\notag\\
&\M(\mathcal{TT}_3)=2^2\cdot 13;\notag\\
&\M(\mathcal{TT}_4)=2^2\cdot5^2\cdot 31;\notag\\
&\M(\mathcal{TT}_5)=2\cdot 3^2\cdot 19^2\cdot 37;\notag\\
&\M(\mathcal{TT}_6)=2\cdot 7^3\cdot 13 \cdot 43 \cdot 127;\notag\\
&\M(\mathcal{TT}_7)=2^7 \cdot 3^5\cdot 5^3\cdot 7\cdot 13\cdot 73.\notag
\end{align}

The amount of factorization is remarkable (the authors have computed and factored
$\M(\mathcal{TT}_n)$ up to $n = 30$) and comparable  to that of the numbers enumerating domino tilings
of squares (given by the well-known formula of Kasteleyn \cite{Kas} and Temperley and Fisher \cite{Fish}). Based on this observation, they posed the following problem.

\begin{prob}[Problem 1.5  in \cite{CK02}]\label{problem24}
 Find a formula for the number of lozenge tilings of $\mathcal{TT}_n$ that explains a
large amount of prime factorization of these numbers.
\end{prob}

\begin{figure}\centering
\includegraphics[width=10cm]{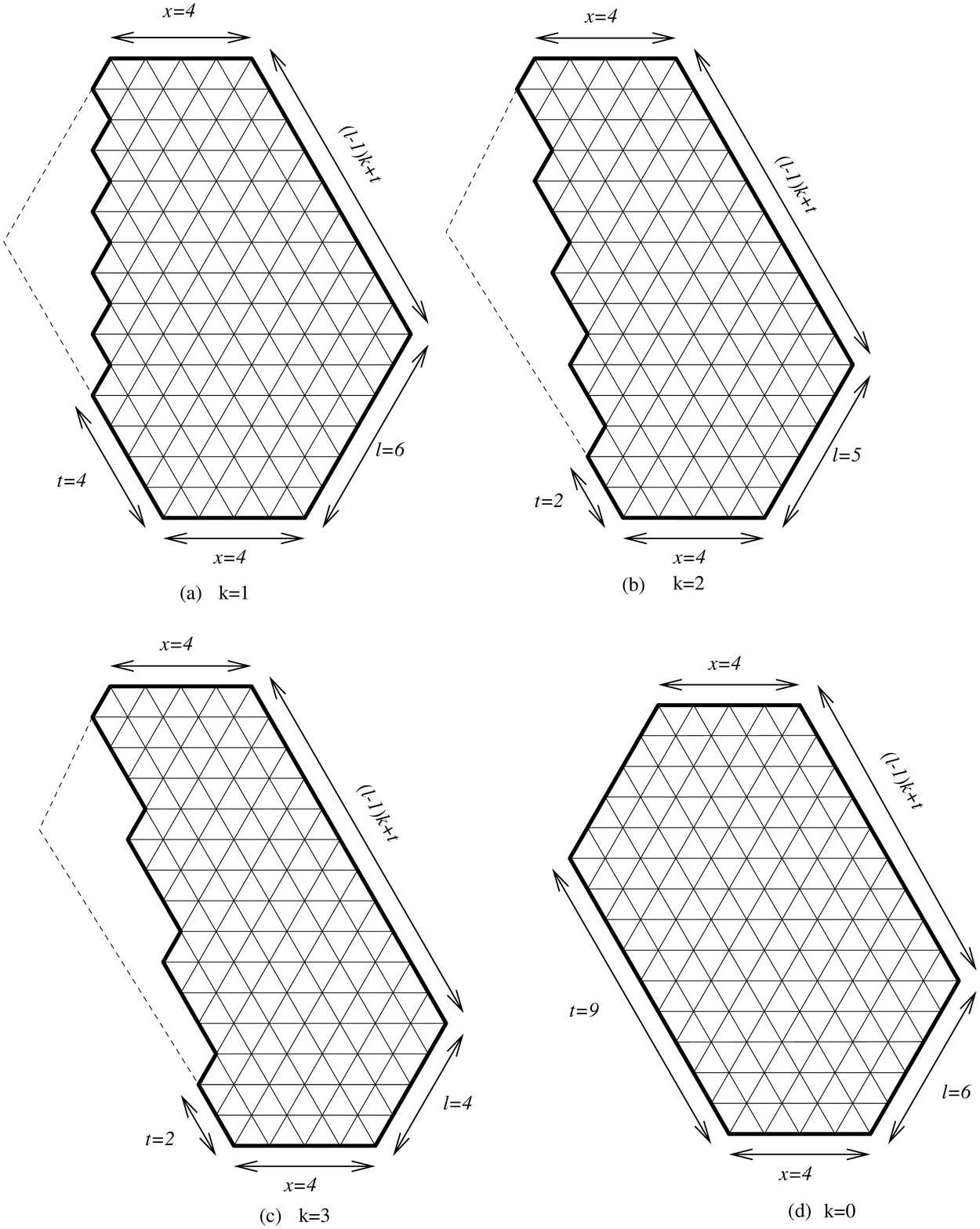}
\caption{The $k$-halved hexagons with no dent.}\label{Fig:Tiltinghex}
\end{figure}

In a recent paper \cite{Tiltinghex}, a common generalization for three famous families of regions in the enumeration of tilings, namely the hexagons, the semi-hexagons, and the halved hexagons, has been introduced. We consider a hexagon, then cut off a maximal $k$-staircase whose each step has the width $k$. When $k=1$, we have exactly the halved hexagon; when $k=0$,   nothing is cut off, and the region is still the hexagon. When $k\geq 2$, we have new regions that are similar to the halved hexagons; however, the cut is tilted. We call the new regions the \emph{$k$-halved hexagons} or the \emph{tilted halved hexagons} (see Figure \ref{Fig:Tiltinghex}).

\begin{figure}\centering
\setlength{\unitlength}{3947sp}%
\begingroup\makeatletter\ifx\SetFigFont\undefined%
\gdef\SetFigFont#1#2#3#4#5{%
  \reset@font\fontsize{#1}{#2pt}%
  \fontfamily{#3}\fontseries{#4}\fontshape{#5}%
  \selectfont}%
\fi\endgroup%
\resizebox{!}{16cm}{
\begin{picture}(0,0)%
\includegraphics{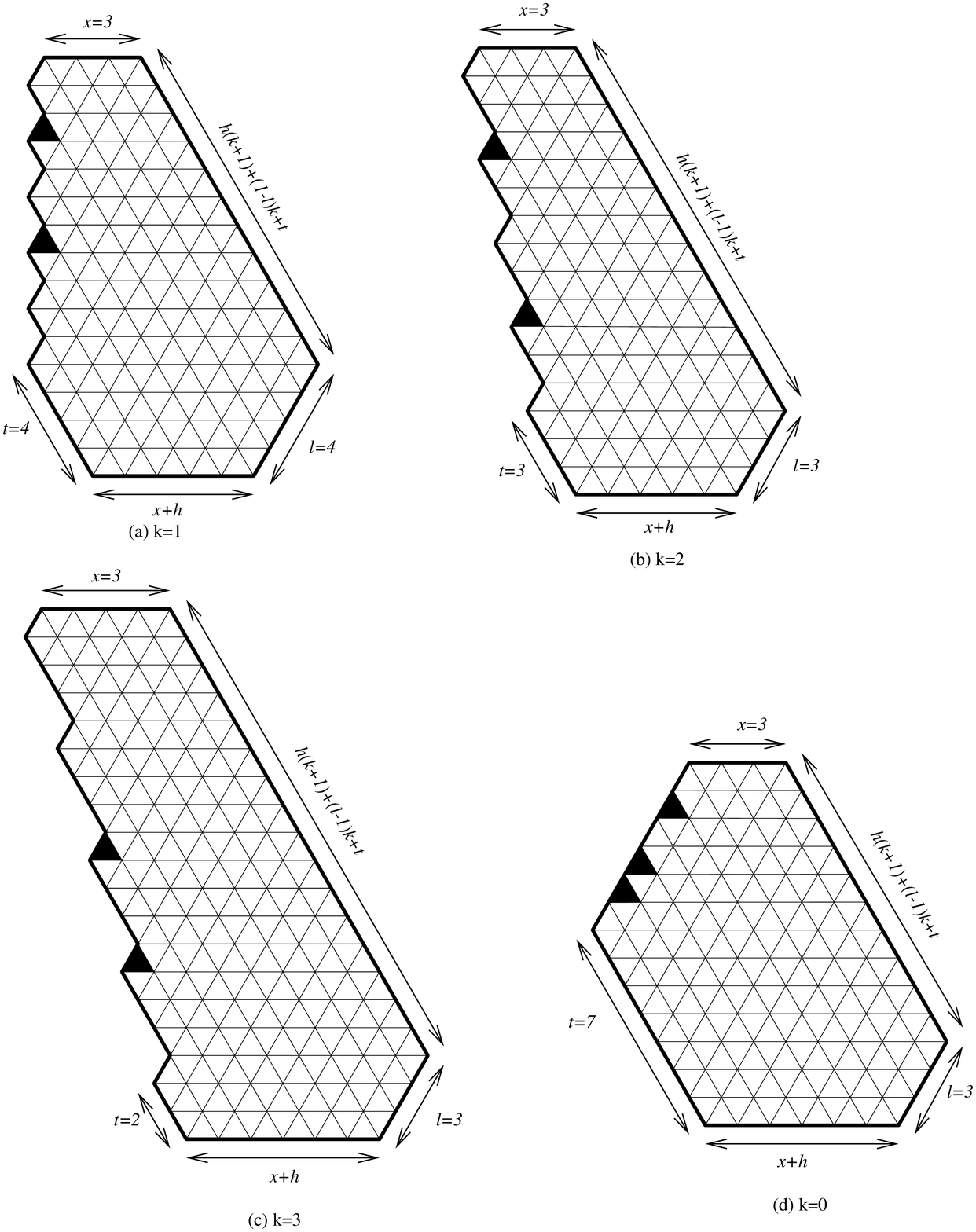}%
\end{picture}%

\begin{picture}(12527,15695)(1225,-15144)
\put(1284,-4179){\makebox(0,0)[lb]{\smash{{\SetFigFont{12}{14.4}{\rmdefault}{\mddefault}{\updefault}{\color[rgb]{0,0,0}$a_1$}%
}}}}
\put(1276,-3466){\makebox(0,0)[lb]{\smash{{\SetFigFont{12}{14.4}{\rmdefault}{\mddefault}{\updefault}{\color[rgb]{0,0,0}$a_2$}%
}}}}
\put(1284,-2041){\makebox(0,0)[lb]{\smash{{\SetFigFont{12}{14.4}{\rmdefault}{\mddefault}{\updefault}{\color[rgb]{0,0,0}$a_3$}%
}}}}
\put(1291,-631){\makebox(0,0)[lb]{\smash{{\SetFigFont{12}{14.4}{\rmdefault}{\mddefault}{\updefault}{\color[rgb]{0,0,0}$a_4$}%
}}}}
\put(7636,-4749){\makebox(0,0)[lb]{\smash{{\SetFigFont{12}{14.4}{\rmdefault}{\mddefault}{\updefault}{\color[rgb]{0,0,0}$a_1$}%
}}}}
\put(2875,-13293){\makebox(0,0)[lb]{\smash{{\SetFigFont{12}{14.4}{\rmdefault}{\mddefault}{\updefault}{\color[rgb]{0,0,0}$a_1$}%
}}}}
\put(8444,-11321){\makebox(0,0)[lb]{\smash{{\SetFigFont{12}{14.4}{\rmdefault}{\mddefault}{\updefault}{\color[rgb]{0,0,0}$a_1$}%
}}}}
\put(7216,-2627){\makebox(0,0)[lb]{\smash{{\SetFigFont{12}{14.4}{\rmdefault}{\mddefault}{\updefault}{\color[rgb]{0,0,0}$a_2$}%
}}}}
\put(9067,-10263){\makebox(0,0)[lb]{\smash{{\SetFigFont{12}{14.4}{\rmdefault}{\mddefault}{\updefault}{\color[rgb]{0,0,0}$a_2$}%
}}}}
\put(1653,-9056){\makebox(0,0)[lb]{\smash{{\SetFigFont{12}{14.4}{\rmdefault}{\mddefault}{\updefault}{\color[rgb]{0,0,0}$a_2$}%
}}}}
\put(9464,-9558){\makebox(0,0)[lb]{\smash{{\SetFigFont{12}{14.4}{\rmdefault}{\mddefault}{\updefault}{\color[rgb]{0,0,0}$a_3$}%
}}}}
\put(6804,-489){\makebox(0,0)[lb]{\smash{{\SetFigFont{12}{14.4}{\rmdefault}{\mddefault}{\updefault}{\color[rgb]{0,0,0}$a_3$}%
}}}}
\put(1240,-7601){\makebox(0,0)[lb]{\smash{{\SetFigFont{12}{14.4}{\rmdefault}{\mddefault}{\updefault}{\color[rgb]{0,0,0}$a_3$}%
}}}}
\end{picture}}
\caption{Several $k$-halved hexagons with dents.}\label{Fig:Tiltinghalvehex2}
\end{figure}

We actually consider  a more general situation  when we allow some ``dents" on the  staircase of the $k$-halved hexagon. Label the staircase levels from the bottom to the top by $1,2,\dots,l+h$. We allow removing the up-pointing unit triangles at $h$ corners. Assume that the remaining steps have labels $a_1,a_2,\dots,a_l$ as they appear from bottom to the top. Denote by $H_{x,t,h}(a_1,a_2,\dots,a_l)$ the resulting region (illustrated in Figure \ref{Fig:Tiltinghalvehex2}). 

As there are nice q-enumerations of the quasi-regular hexagons, semi-hexagons, and halved hexagons, we would expect a nice $q$-enumeration of the $k$-halved hexagons.
\begin{prob}\label{problem25}
Find a $q$-enumeration of tilings of a $k$-halved hexagon.
\end{prob}

\begin{figure}\centering
\includegraphics[width=15cm]{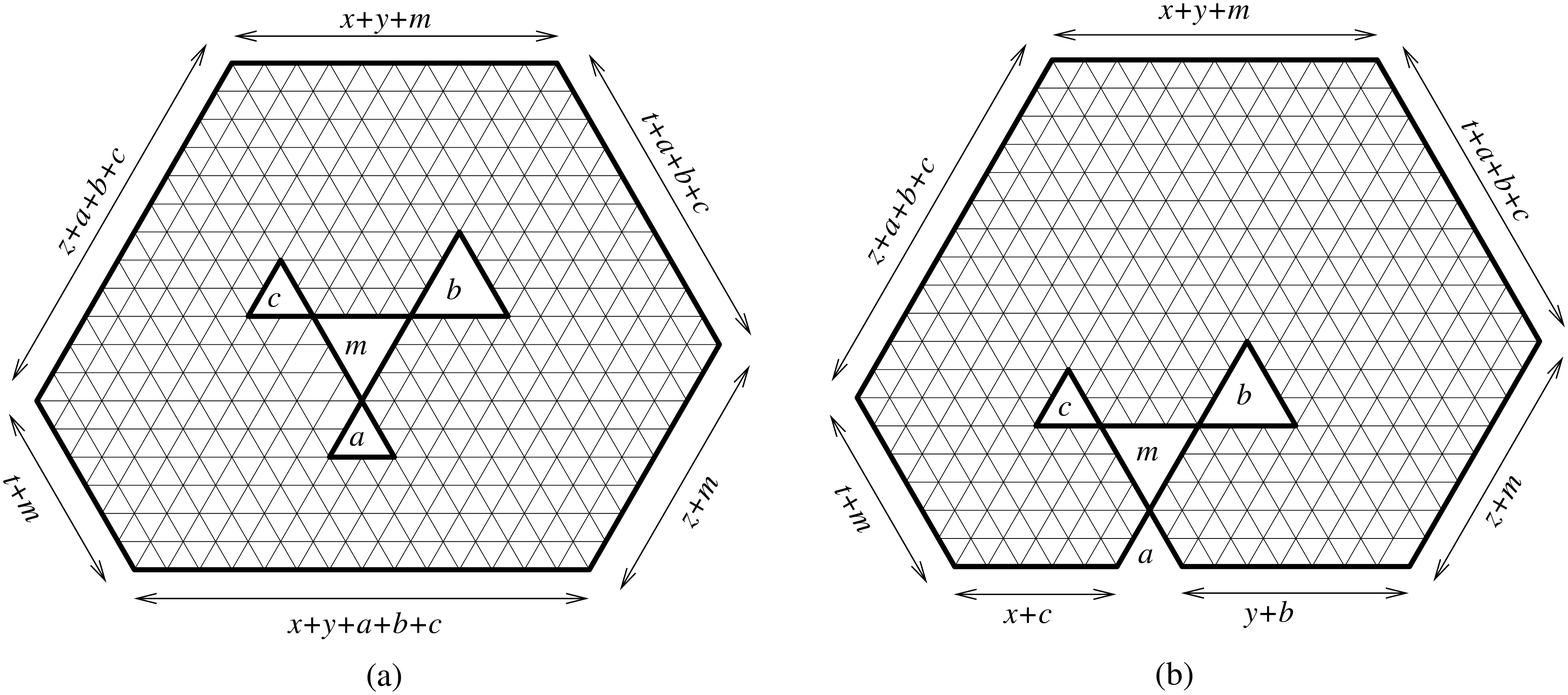}
\caption{(a) A hexagon with a shamrock hole in the center in \cite{CK13}.  (b) A hexagon with a shamrock hole on the boundary in \cite{Tri19}. The picture first appeared in \cite{LR}.}\label{Fig:Shamrock}
\end{figure}

\begin{figure}\centering
\includegraphics[width=12cm]{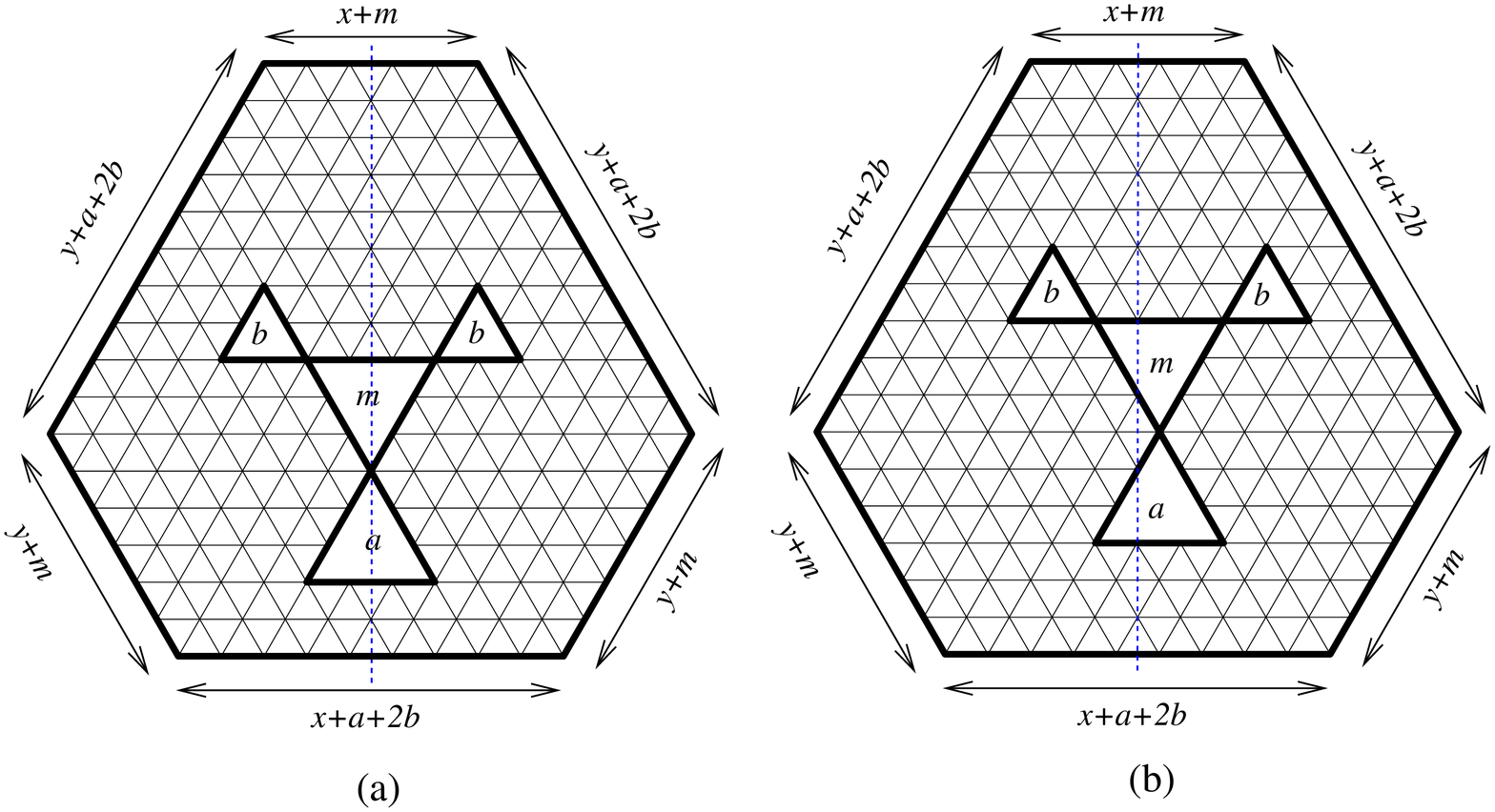}
\caption{(a) Placing a shamrock hole on the symmetry axis of the hexagon. (b) Placing a shamrock hole $1/2$ unit off the symmetry axis of the hexagon.}\label{Fig:offcenter}
\end{figure}

In the enumeration of tilings, the symmetric regions often behave better than the asymmetric ones. Many families of regions may not have a nice tiling number in the general case. However, their tilings are enumerated by simple products in the symmetric case. 
Let us revisit Ciucu--Krattenthaler's $S$-cored hexagon, a hexagon with a cluster of four triangles (called a ``shamrock") removed. The $S$-cored hexagons have a nice tiling number in only two situations: (1) the ``shamrock"  is removed from the center \cite{CK13}, and (2) the shamrock is removed from the boundary \cite{Tri19} (see Figure \ref{Fig:Shamrock} for examples). In general, if we remove the shamrock from a place that is different from the center or the boundary, then the number of tilings is not given by a simple product formula. However, in the case of reflectively symmetric hexagons, it has been shown \cite{LR} that we can remove the shamrock at \emph{any} position along the symmetry axis and still get a beautiful tiling formula (see Figure \ref{Fig:offcenter}(a)). 

\begin{figure}\centering
\includegraphics[width=12cm]{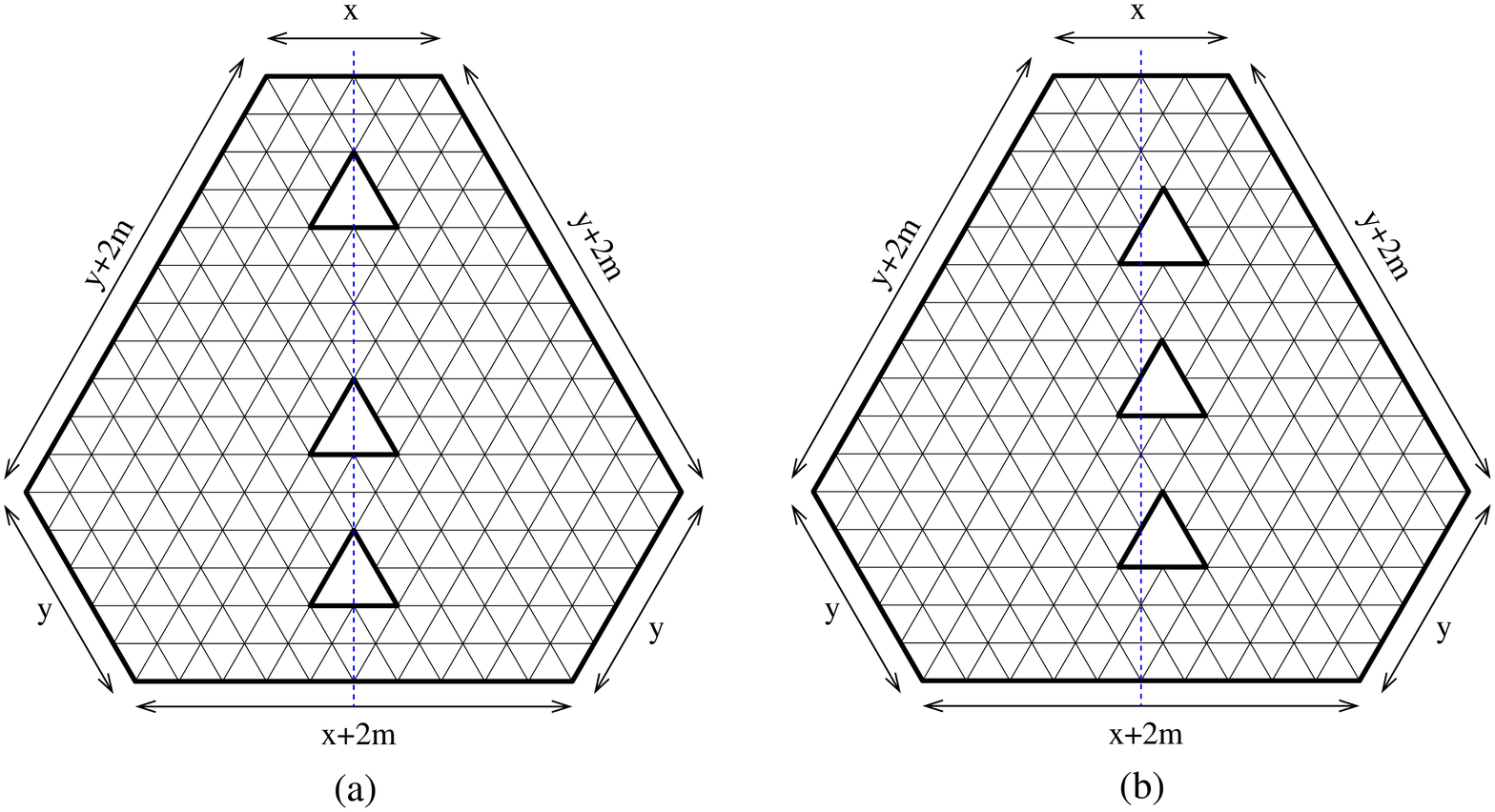}
\caption{(a) Placing triangular holes of side-length 2 on the symmetry axis of the hexagon. (b) Placing triangular holes $1/2$ unit off the symmetry axis.}\label{Fig:offcenter2}
\end{figure}

We have the same observation for hexagons with a family of vertically aligned triangular holes of side-length 2.  In the general case, we do not have a nice tiling number; however, in the case of reflectively symmetric hexagons, Ciucu proved a simple product formula for the tiling number \cite{Ciucu2}. See Figure \ref{Fig:offcenter2}(a) for an illustrated picture. 

 Our data suggests that if we place the shamrock holes and the family of triangular holes $1/2$ unit off the symmetry axis, we still have nice tiling numbers in the above two families of regions. See Figures \ref{Fig:offcenter}(b) and \ref{Fig:offcenter2}(b). In general, we often have the following ``\emph{off-center phenomenon}": if we place the hole(s) not on the symmetry axis but $1/2$ unit off the symmetry axis, then the number of tilings seems to be as nice as the tiling number in the symmetric case. This fact is not true anymore if we move the holes just 1 unit away from the symmetric axis.  It would be interesting to have a precise explanation for this phenomenon.

\begin{prob}\label{problem26}
Explain the off-center phenomenon.
\end{prob}

It is worth noticing that the enumeration of tilings of  regions with holes\footnote{Several authors used the term ``\emph{holey regions}" for  ``regions with holes".} is especially challenging. One of the difficulties is that the  Lindstr\"{o}m--Gessel--Viennot determinant  does \emph{not} give the tiling number; it gives the \emph{signed} tiling number instead. We prefer the reader to, say, the proofs in Section 4 of \cite{CEKZ} for a detailed explanation. In the reflectively symmetric case, one can go around this obstacle by using a powerful tool, usually mentioned as \emph{Ciucu's factorization theorem} (see \cite[Theorem 1.2]{Ciucu97}). This theorem allows us to simplify the case of symmetric regions with holes to the case of simply connected regions, say by dividing the region along the symmetric axis into two smaller regions with no holes. See \cite[Section 3]{Ciucu2} for more details of the method. However, if we slide the holes $1/2$ unit away from the symmetry axis, then Ciucu's method is \emph{failed} to apply (as the new region is not symmetric anymore). It would be very interesting to find explicit formulas for the numbers of tilings of  the two regions in Figures \ref{Fig:offcenter}(b) and \ref{Fig:offcenter2}(b).

\begin{prob}\label{problem27}
 Find a formula for the number of tilings of a symmetric hexagon with a shamrock hole at $1/2$ unit off the symmetric axis (as in Figure \ref{Fig:offcenter}(b)).
\end{prob}

\begin{prob}\label{problem28}
 Find a formula for the number of tilings of a symmetric hexagon with a family of aligned $2$-triangles at $1/2$ unit off the symmetric axis (as in Figure \ref{Fig:offcenter2}(b)).
\end{prob}

\begin{figure}\centering
\includegraphics[width=15cm]{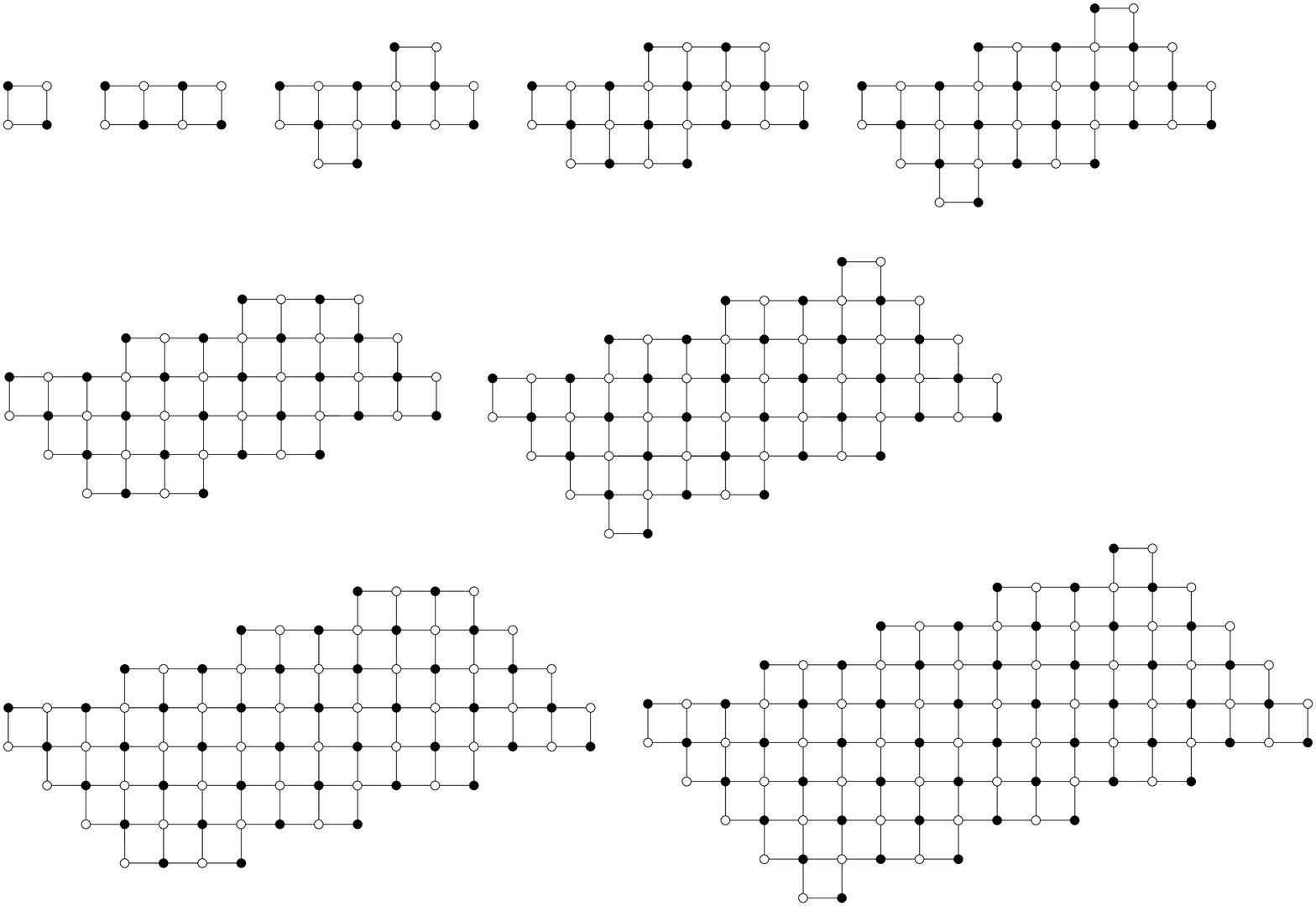}
\caption{The Aztec pillows $AP_n$, for $n=1,2,3,\dots,9$ (reading from left to right, top to bottom).}\label{Fig:Aztecpillow}
\end{figure}

The ``\emph{Aztec pillow graphs}" were first introduced in \cite{Propp}. The width of an Aztec pillow is always even. When the width of an Aztec pillow is $4k+2$, the upper half of the graph has $k+1$ up-steps ($k$ steps of size $3$, and one step of size $1$), followed by $k+1$ down-steps of size 1 (as we go from the left to right), and the lower half is simply a $180^{\circ}$-rotation of the upper half. When the Aztec pillow's width is $4k$, the upper half has $k$ up-steps of size 3, followed by $k$ down-steps of size $1$. Denote by $AP_n$ the Aztec pillow of width $2n$. See Figure \ref{Fig:Aztecpillow} for the graphs $AP_n$, for $n=1,2,\dots,9$. Forest Tong \cite{Tong} conjectured the following elegant divisibility of the matching numbers of the Aztec pillows.

\begin{prob}\label{problem29}
Prove that $\MM(AP_m)\, |\, \MM(AP_n)$ whenever $(m + 3)\, | \, (n + 3)$, where $\MM(G)$ denotes the number of perfect matching of graph $G$.
\end{prob}
Tong  has verified this conjecture computationally for $m, n < 77$. 

It is worth noticing that the Aztec pillows' definition has been generalized by Christopher Hanusa. The original Aztec pillows above are the \emph{Aztec 3-pillows}. Forest Tong also observed that \emph{none} of the 5-pillows, 7-pillows, and 9-pillows (as defined in \cite{Hanusa}) seem to share the above divisibility property.

\medskip

Next, we investigate a nice property of the \emph{matching polynomial} (see, e.g., \cite[p.  xxxii,p. 333ff]{LoPl}). 

Let $G$ be a graph with no loops. Denote by $m_k(G)$ the number of (partial) matchings of the graph $G$ with exactly $k$ edges, where $m_0(G)=1$ by convention. The \emph{matching polynomial} of $G$ is defined to be
\[\mathcal{M}(G)=\sum_{k\geq 0}m_k(G) x^k.\]

A sequence $(a_i)_{i\geq 0}$ is a  \emph{P\'{o}lya frequency sequence} (PFS) if  the infinite Toeplitz matrix $(M_{i,j})_{i,j\geq 0}$ defined by $M_{i,j}=a_{j-i}$ (where $a_k=0$ if $k<0$ by convention), i.e,
\[(M_{i,j})_{i,j\geq 0}=\begin{pmatrix} a_0 &a_1 &a_2&a_3&a_4&\cdots\\
0&a_0 &a_1 &a_2&a_3&\cdots\\
0&0&a_0 &a_1 &a_2&\cdots\\
0&0&0&a_0 &a_1 &\cdots\\
0&0&0&0&a_0 &\cdots\\
0&0&0&0&0 &\cdots\\
\vdots&\vdots&\vdots &\vdots &\vdots&\ddots
\end{pmatrix},\]
 has all nonnegative minors. By definition, a P\'{o}lya frequency sequence is log-concave, as the log-concavity is equivalent to the fact that all $2\times 2$ minors of $(M_{i,j})_{i,j\geq 0}$ are nonnegative.

O. J. Heilmann and  E. H. Lieb  \cite[Theorem  4.2]{HL} proved that the matching polynomial $\mathcal{M}(G)$ has all real roots. This implies that the sequence $(m_k(G))_{k\geq 0}$ is a P\'{o}lya frequency sequence (see, e.g., \cite{Brenti}).    It would be  interesting to find a  combinatorial proof for this fact.

\begin{prob}\label{problem30}
Prove combinatorially that the sequence of matching numbers $(m_k(G))_{k\geq 0}$ is a P\'{o}lya frequency sequence.\end{prob}
We note that Krattenthaler \cite{Kra96} provided a combinatorial proof for a special case of Problem \ref{problem30}, namely the log-concavity of the sequence $(m_k(G))_{k\geq 0}$.

We can generalize Problem \ref{problem30} to weighted graphs as follows. To each edge $e$ of $G$, we assign a weight $x_e$. Then we define the weighted matching number $m_k(G,\textbf{x})$ by
\[m_k(G,\textbf{x})=\sum_{M}\prod_{e\in M} x_e,\]
where the sum is taken over all $k$-element matchings $M$ in $G$. Assume $P_i=P_i(\textbf{x})$ be a polynomial in $\mathbb{Z}[\textbf{x}]$, for $i=1,2,3,\dots$. We now define a sequence $(P_i)_{i=1}^{\infty}$ to be an \emph{$\textbf{x}$-PFS} if all minors of the Toeplitz matrix
$(P_{j-i})_{i,j\geq 0}$ (where $P_{k}=0$ if $n<0$) are polynomials in the $x_e$'s with nonnegative coefficients. 
\begin{prob}\label{problem31}
Prove combinatorially that the sequence of $\textbf{x}$-matching numbers $(m_k(G,\textbf{x}))_{k\geq 0}$ is an $\textbf{x}$-PFS.
\end{prob}

Recently, P. Galashin and P. Pylyavskyy \cite{GP} consider a similar positivity for a planar bipartite graph $G=(V_1,V_2,E)$ embedded on a cylinder $\mathcal{O}$. Let $\tau_1$ and $\tau_2$ be two perfect matchings on $G$. We always orient the edges in a perfect matching $\tau$ of $G$ from a vertex in $V_1$ to a vertex in $V_2$. Define the difference of two perfect matching $\tau_1-\tau_2$ to be the directed graph on $\mathcal{O}$ with vertices $V_1\cup V_2$ obtained by superimposing $\tau_1$ and $\tau_2$, and the reversed the direction of edges in $\tau_2$.  This way $\tau_1-\tau_2$ is always a disjoint union of directed simple cycles, which can be viewed as singular 1-cycles on $\mathcal{O}$. We define the \emph{relative height} of two perfect matchings $h(\tau_1,\tau_2)$ to be the image of $H_1(\mathcal{O},\mathbb{Z})\simeq \mathbb{Z}$ of the sum of these cycles. Fix a minimal-height perfect matching $\tau_0$. Next, we can define the \emph{absolute height} of a perfect matching as $h(\tau):=h(\tau,\tau_0)$.

We also assume that the edges of $G$ are weighted by $x_e$'s as above. We now define $H_i(\textbf{x})$ to be the sum of weights of all perfect matchings with height $i$, where the weight of  a perfect matching is the product of its edge-weights  as usual.

\begin{prob}[Conjecture 6.1 in \cite{GP}]\label{problem32}
Prove that the sequence $(H_i(\textbf{x}))_{i\geq 0}$ is an $\textbf{x}$-PFS.
\end{prob}
We note that while the $PFS$ properties in Problems \ref{problem30} and \ref{problem31} have been proved (and we are asking for a combinatorial proof), the question about the PFS property in Problem \ref{problem32} is still open. Several special cases of Problem \ref{problem32} have been proved in \cite{CPS}.

\medskip

We conclude this section with a problem of a rather different flavor. Let $k$ be a fixed positive integer. Denote by $A_n=A_{n,k}$ the number of domino tilings of a $k\times n$ rectangle. Form the generating function
\[F_k(x)=\sum_{n\geq 0}A_n x^n.\]
It has been shown that $F_k(x)$ can be written as a rational function, say $F_k(x)=\frac{P_k(x)}{Q_k(x)}$ with $P_k$ and $Q_k$ polynomials with integer coefficients, and $Q_k(0)=1$ \cite{Klarner2}. Stanley proves that all the roots of $Q_k(x)$ are real and nonzero, and exactly half of the roots are positive \cite{Stanley4}. He also conjectures a special pattern for the roots of the polynomial:

\begin{prob}\label{problem35}
Prove that $Q_k(x)$ has distinct roots.
\end{prob}

\section*{Acknowledgement} The author would like to thank Gregg Musiker, James Propp,  Pavlo Pylyavskyy, Richard Stanley, and Dennis Staton for helpful comments and fruitful discussions. Problem \ref{problem6} was suggested by Dennis Stanton, and Problems \ref{problem30}--\ref{problem32} were introduced to the author by Pavlo Pylyavskyy. 

\bibliographystyle{plain}
\bibliography{TilingOPAC2022}

\end{document}